\documentclass[11pt]{article}
\usepackage{amsmath,amssymb,amsthm}
\allowdisplaybreaks[4]
\usepackage{xcolor}
\usepackage{graphicx}
\usepackage{tikz}
\usepackage{float}
\usepackage[T1]{fontenc}
\usepackage[a4paper,margin=1in]{geometry}
\usepackage[english]{babel}
\usepackage{newtxtext,newtxmath}
\usepackage{hyperref}
\hypersetup{
	hidelinks,
	pdftitle={Modular Ideals and Factorizations of Affine Hyperplane Arrangements},
	pdfauthor={Yanru Chen, Weikang Liang, Suijie Wang, and Chengdong Zhao},
	pdfsubject={Hyperplane arrangements and intersection semilattices},
	pdfkeywords={hyperplane arrangements, intersection semilattices, ideal decompositions, modular ideals, characteristic polynomials, Orlik-Solomon algebras, geometric realizations, semimatroids}
}

\tikzset{
	hasse node/.style={draw, rounded corners=1.2pt, inner sep=1.5pt,
		minimum height=4.8mm, minimum width=6mm, font=\scriptsize},
	hasse blue/.style={hasse node, fill=blue!16, draw=blue!70!black},
	hasse red/.style={hasse node, fill=red!16, draw=red!70!black},
	hasse green/.style={hasse node, fill=green!16, draw=green!55!black},
	hasse gray/.style={hasse node, fill=gray!12, draw=gray!65!black},
	hasse bottom/.style={hasse node, fill=gray!8, draw=gray!65!black},
	hasse cover/.style={draw=gray!58, line width=0.36pt},
	hasse blue cover/.style={draw=blue!65!black, line width=0.52pt},
	hasse red cover/.style={draw=red!65!black, line width=0.52pt},
	hasse green cover/.style={draw=green!55!black, line width=0.52pt}
}

\newtheorem{theorem}{Theorem}[section]
\newtheorem{lemma}[theorem]{Lemma}
\newtheorem{prop}[theorem]{Proposition}
\newtheorem{coro}[theorem]{Corollary}
\newtheorem{fact}[theorem]{Fact}
\theoremstyle{definition}
\newtheorem{defi}[theorem]{Definition}
\theoremstyle{remark}
\newtheorem*{remark*}{Remark}
\newtheorem{example}[theorem]{Example}

\newcommand{\A}{\mathcal{A}}
\newcommand{\B}{\mathcal{B}}
\numberwithin{equation}{section}
\begin{document}

\begin{center}
\setlength{\parskip}{0pt}
{\Large\bfseries
Modular Ideals and Factorizations of
Affine Hyperplane Arrangements\par}

\vspace{7pt}

Yanru Chen$^{1}$,
Weikang Liang$^{2,*}$,
Suijie Wang$^{3}$
and Chengdong Zhao$^{4}$\par

\vspace{8pt}

$^{1,2,3}$ School of Mathematics\\
Hunan University\\
Changsha 410082, Hunan, P. R. China\par

\vspace{12pt}

$^{4}$ School of Mathematics and Statistics\\
Central South University\\
Changsha 410083, Hunan, P. R. China\par

\vspace{15pt}

\begin{tabular}{@{}r@{\ }l@{}}
Emails: &
$^{1}$\href{mailto:yanruchen@hnu.edu.cn}{yanruchen@hnu.edu.cn},
$^{2}$\href{mailto:kangkang@hnu.edu.cn}{kangkang@hnu.edu.cn},\\
&
$^{3}$\href{mailto:wangsuijie@hnu.edu.cn}{wangsuijie@hnu.edu.cn},
$^{4}$\href{mailto:cdzhao@csu.edu.cn}{cdzhao@csu.edu.cn}
\end{tabular}\\[2pt]
$^{*}$ Corresponding author\par
\end{center}

\vspace{3mm}
\begin{abstract}
Let \(\A\) be an affine hyperplane arrangement, let \(L(\A)\) be its
intersection semilattice, and let \(\chi_{\A}(t)\) be its characteristic
polynomial. We introduce ideal decompositions of finite ranked
meet-semilattices and prove that they induce purely combinatorial
factorizations of characteristic polynomials. This framework recovers
both Stanley's factorization associated with modular elements and Terao's
factorization associated with nice partitions. For simple semimatroids,
every join-closed ideal decomposition comes from a direct-sum
decomposition. The two-factor case defines modular ideals, extending the
role of modular elements from central to affine arrangements.

Over an infinite field, every modular ideal of an arrangement intersection
semilattice is realized by a suitable translated restriction of a
subarrangement. We exhibit a nonempty Zariski-open set of valid
translations, yielding essential realizations of the given modular ideal.
We further prove that modular subarrangements correspond under
coning to modular elements on the hyperplane at infinity and hence to
M-ideals in the affine setting. Finally, every modular ideal gives an
Orlik--Solomon graded vector-space decomposition, which becomes a
graded-algebra decomposition when the two factors come from
subarrangements.

\end{abstract}

\noindent\textbf{Keywords.}\quad
hyperplane arrangement;
intersection semilattice;
ideal decomposition;
modular ideal;
characteristic polynomial;
Orlik--Solomon algebra;
geometric realization;
semimatroid.

\medskip
\noindent\textbf{2020 Mathematics Subject Classification.}\quad
Primary 52C35; Secondary 05B35, 06A07, 16W50.

\section{Introduction and main results}

An affine \emph{hyperplane arrangement} \(\A\) is a finite collection of
affine hyperplanes in \(V\cong\mathbb{K}^d\) over a field \(\mathbb{K}\).
Its \emph{intersection semilattice} (also called the intersection poset)
\(L(\A)\) consists of all nonempty intersections of hyperplanes in
\(\A\), ordered by reverse inclusion: \(X\leq Y\) if and only if
\(Y\subseteq X\). The ambient space \(V\) is the minimal element,
corresponding to the empty intersection. The poset \(L(\A)\) is a finite
ranked meet-semilattice with rank function
\(r(X)=d-\dim(X)\).

The arrangement \(\A\) is \emph{central} if the intersection of all
hyperplanes in \(\A\) is nonempty; then \(L(\A)\) is a lattice. Throughout,
\(\A\) is assumed to be \emph{essential}, meaning that \(L(\A)\) has rank
\(d\). Accordingly, we use that \(L(\A)\) is graded, so every maximal
element has rank \(d\); see \cite{Orlik-Terao,Wachs}. The
\emph{characteristic polynomial} of \(\A\) is
\[
\chi_\A(t) = \sum_{X \in L(\A)} \mu(X) t^{\dim(X)},
\]
where $\mu$ is the M\"{o}bius function on $L(\A)$.
For standard background on intersection semilattices and characteristic
polynomials of arrangements, see \cite{Orlik-Terao,Stanley-book}.

Our starting point is an \emph{ideal decomposition}: a collection of
order ideals satisfying a local principal-ideal condition, such that
joins of elements from different ideals exist and their ranks add. This
notion is defined in Section
\ref{sec:ideal-decompositions}. It yields a M\"obius convolution and a
characteristic-polynomial factorization for finite ranked
meet-semilattices. For an ideal \(I\), define its \emph{orthogonal ideal} by
\[
I^\perp:=\{y\in L:y\wedge x=\hat{0}\text{ for every }x\in I\}.
\]
We call \(I\) a \emph{modular ideal} when \(\{I,I^\perp\}\) is an ideal
decomposition. We prove that modular ideals admit translated geometric
realizations; they also yield, via multiplication, decompositions of
Orlik--Solomon algebras. In such a decomposition, the factor indexed by
an ideal is a subalgebra exactly when that ideal is join-closed.

\paragraph{Relation to previous work.}
For geometric lattices, Stanley's modular-element theorem factors the
characteristic polynomial \cite{Stanley1971}; Brylawski related this
factorization to complete Brown truncation \cite{Brylawski}, and Liu and
Sagan studied left-modular extensions \cite{LiuSagan}. Terao connected
modular flats with arrangement fibrations \cite{HT1986,Orlik-Terao};
fiber-type arrangements were further developed by Falk--Randell and
Jambu \cite{FalkRandell,Jambu}. These results apply directly to central
arrangements. Proposition \ref{prop:central-modular-ideal} recovers
Stanley's theorem, while modular ideals extend the factorization to
affine intersection semilattices.

Terao's nice partitions give graded vector-space factorizations of
Orlik--Solomon algebras and linear factorizations of characteristic
polynomials \cite{HT1992}; addition--deletion and inductive versions
appear in \cite{JambuParis,HogeRoehrle}. Proposition
\ref{prop:nice-partition-ideal-decomposition} identifies nice partitions
with ideal decompositions generated by hyperplane blocks. Hallam and
Sagan's quotient-poset method instead partitions the atoms and passes to
a quotient poset \cite{Stanley1972,Sagan15}. In our setting, the
factors remain order ideals of the original semilattice.

Geometric semilattices and semimatroids provide the natural
combinatorial setting for affine arrangements
\cite{Wachs,Ardila-semimatroids}; see also Kawahara's quasimatroid
framework \cite{Kawahara}. Theorem
\ref{thm:join-closed-semimatroid-direct-sum} shows that a join-closed
ideal decomposition of a simple semimatroid comes from a direct sum. The
term \emph{modular ideal} used here refers to an order ideal of the flat
semilattice, not to the set-system notion of a modular ideal of a matroid
\cite[Definition~3.3]{Ardila-semimatroids}.

Finally, Bibby and Delucchi's M-ideals extend supersolvability in locally
geometric posets \cite{Bibby}; related settings appear in
\cite{Bibby-1,Bibby-matroid-schemes}. In the affine
arrangement setting, M-ideals correspond to the join-closed modular
ideals induced by subarrangements. General modular ideals need not be
join-closed, as Example \ref{ex:modular-ideal-realization} shows.

For a finite ranked meet-semilattice \(P\), and also for an ideal
equipped with the restricted rank, write
\[
r(P):=\max_{x\in P}r(x),
\qquad
\chi(P,t):=\sum_{x\in P}\mu(\hat{0},x)t^{r(P)-r(x)}.
\]

\begin{theorem}\label{thm:ideal-factorization}
Let $L$ be a finite ranked meet-semilattice.
If $L$ admits an ideal decomposition $\mathcal{I} = \{I_1, I_2, \ldots, I_m\}$, then
\[
\chi(L,t) = t^{c} \prod_{i=1}^{m}\chi(I_{i},t),
\]
where $c = r(L) - \sum_{i=1}^{m} r(I_i)$.
\end{theorem}

Theorem \ref{thm:ideal-factorization} is the combinatorial core of the
paper. It requires neither a representation nor an actual product
decomposition of \(L\). For intersection semilattices, the rank
correction vanishes; nice partitions are recovered as the special case
generated by hyperplane blocks (Proposition
\ref{prop:nice-partition-ideal-decomposition}).

For an ideal \(I\subseteq L(\A)\), identify each hyperplane with its
rank-one flat in \(L(\A)\), and set
\[
\A_I:=\{H\in\A\mid H\in I\}.
\]
The next result realizes every modular ideal by restricting \(\A_I\) to
a suitable translate of a maximal element of \(I^\perp\).
For an arrangement \(\B\) in \(V\) and an affine subspace
\(Y\subseteq V\), the \emph{restriction} of \(\B\) to \(Y\) is
\begin{align*}
	\B^{Y} = \{H \cap Y \mid H \in \B,\;H \cap Y \neq \emptyset,\;H\cap Y\neq Y\}.
\end{align*}
For \(\alpha\in V\), write
\[
Y + \alpha = \{y + \alpha \mid y \in Y\},
\]
which is again an affine subspace of \(V\).
\begin{theorem}\label{thm:geometric-realization}
Let $\A$ be an affine hyperplane arrangement over an infinite field and let
$L(\A)$ be its intersection semilattice. If $I$ is a modular ideal of
$L(\A)$, then for every maximal element $Y$ of $I^\perp$, there exists
$\alpha \in V$ such that
\[
I \cong L\left((\A_I)^{Y + \alpha}\right).
\]
\end{theorem}

Theorem \ref{thm:realization-parameter-space} describes the valid
translations as a nonempty Zariski-open set and gives a finite-field
criterion. For a central arrangement, the same construction realizes the
quotient in Stanley's modular factorization; see Corollary
\ref{cor:central-restriction-realization}.

A modular ideal need not equal \(L(\A_I)\). This motivates
\emph{modular subarrangements}, whose intersection semilattices are
modular ideals. Theorem \ref{thm:modular-subarrangements-cone}, together
with \cite{Bibby}, gives
\[
\begin{aligned}
\{\text{M-ideals of }L(\A)\}
&\longleftrightarrow
\{\text{modular subarrangements of }\A\}\\
&\subseteq
\{\text{modular ideals of }L(\A)\}.
\end{aligned}
\]
The inclusion can be strict (Example
\ref{ex:modular-ideal-realization}). Thus subarrangements and M-ideals
capture only a subclass of modular ideals; characteristic-polynomial
factorization and translated realization hold on the larger class.

The final main result is the affine analogue of Terao's multiplication
theorem \cite{HT1986}. Its conclusion is linear in general: the two
factors need not be subalgebras.
We use the standard notation \(A_X(\A)\) for the local Brieskorn
component indexed by \(X\in L(\A)\); the definition is recalled at the
beginning of Section~\ref{sec:os-factorization}.

\begin{theorem}\label{thm:os-factorization}
	Let \(J\) be a modular ideal of \(L(\A)\), and take the
	Orlik--Solomon algebra over a coefficient field \(\Bbbk\). Then the
	graded \(\Bbbk\)-linear map
	\begin{align}\label{eq:os-multiplication-map}
		\kappa \colon
		\Big(\bigoplus_{X \in J}A_X(\A)\Big) \otimes
		\Big(\bigoplus_{Y \in J^\perp}A_Y(\A)\Big)
		\longrightarrow A(\A)
	\end{align}
	defined by multiplication in \(A(\A)\) is an isomorphism of graded
	\(\Bbbk\)-vector spaces.
\end{theorem}

Theorem \ref{thm:subalgebra-criterion} shows that the factor indexed by
\(J\) is a subalgebra exactly when \(J=L(\A_J)\). In that case
\(A(\A)\) is free over \(A(\A_J)\); if the same holds for \(J^\perp\),
multiplication is an isomorphism of graded algebras.

The paper is organized as follows. Section
\ref{sec:ideal-decompositions} proves the M\"obius convolution and its
characteristic-polynomial consequences, including the nice-partition and
semimatroid cases. Section \ref{sec:modular-realization} develops modular
ideals and their translated realizations. Section
\ref{sec:modular-subarrangements} treats modular subarrangements and
M-ideals. Section
\ref{sec:os-factorization} proves the Orlik--Solomon decomposition and
the subalgebra criterion.

\section{Ideal decompositions and M\"obius convolution}
\label{sec:ideal-decompositions}

We first prove the M\"obius convolution underlying Theorem
\ref{thm:ideal-factorization}, then treat intersection semilattices, nice
partitions, and join-closed semimatroid decompositions.

\subsection{Ideal decompositions of ranked meet-semilattices}
\label{subsec:ranked-meet-semilattices}
We use standard terminology for finite posets; see
\cite{Stanley-book}. For \(x,y\) in a finite poset \(P\), write
\(x\lessdot y\) if \(y\) covers \(x\). When they exist, the join and meet
are denoted by \(x\vee y\) and \(x\wedge y\), respectively. A
\emph{meet-semilattice} is a poset in which every pair has a meet, and a
\emph{lattice} is a poset in which every pair has both a join and a meet.
Thus a finite meet-semilattice has a unique minimal element \(\hat{0}\);
a finite lattice also has a unique maximal element \(\hat{1}\).

A finite meet-semilattice \(L\) is \emph{ranked} if it has a rank function
\(r\colon L\to \mathbb{N}\) satisfying \(r(\hat{0})=0\) and
\(r(y)=r(x)+1\) whenever \(x\lessdot y\). Its rank is
\(r(L):=\max\{r(x)\mid x\in L\}\). It is \emph{graded} if all maximal
elements have rank \(r(L)\).

For a finite ranked meet-semilattice \(L\),
the \textit{M\"{o}bius function} $\mu \colon L \times L \longrightarrow \mathbb{Z}$ is defined recursively by
\[
\mu(x, y) =
\begin{cases}
1, & \text{if } x = y, \\[3pt]
-\sum\limits_{x \leq z < y} \mu(x, z), & \text{if } x < y, \\[3pt]
0, & \text{otherwise}.
\end{cases}
\]
In particular, for any \(x\in L\), we write
$
\mu(x) = \mu(\hat{0}, x)
$
for brevity.
The \emph{characteristic polynomial} of \(L\) is defined by
\begin{align}\label{eq:poset-characteristic-polynomial}
\chi(L, t) = \sum_{x \in L}\mu(x)t^{r(L) - r(x)}.
\end{align}

Throughout this section, \(L\) is a finite ranked meet-semilattice with
rank function \(r\). An \emph{ideal} of \(L\) is a nonempty downward closed subset:
if \(x\in I\) and \(y\leq x\), then \(y\in I\). For
\(A=\{x_1,\ldots,x_k\}\subseteq L\), let
\[
\langle A\rangle=\langle x_1,\ldots,x_k\rangle
:= \{y\in L\mid y\leq x_i\text{ for some }i\}.
\]
A \emph{principal ideal} is one generated by a singleton; the ideal
\(\{\hat{0}\}\) is called trivial. Any ideal inherits the restricted rank
function, again denoted by \(r\) when no confusion can arise. Moreover,
for \(x\) in an ideal \(I\), the interval \([\hat{0},x]\) is the same in
\(I\) and in \(L\); hence the corresponding M\"obius values agree.

\begin{defi}\label{def:ideal-decomposition}
	A collection \(\mathcal{I}=\{I_1,I_2,\ldots,I_m\}\) of ideals of \(L\)
	is called an \emph{ideal decomposition} of \(L\) if the following three
	conditions hold:
	\begin{itemize}
		\item[\hypertarget{meet_condition}{$(\mathrm{I1})$}] \textbf{(meet condition)} For every \(x\in L\setminus\{\hat{0}\}\), there is an index \(j\in\{1,\ldots,m\}\) such that
		$I_j \cap \langle x \rangle$ is a nontrivial principal ideal of $L$;
		\item[\hypertarget{join_condition}{$(\mathrm{I2})$}] \textbf{(join condition)} For every choice \(x_i\in I_i\), \(1\leq i\leq m\), the join \(x_1\vee\cdots\vee x_m\) exists;
		\item[\hypertarget{rank_condition}{$(\mathrm{I3})$}] {\textbf{(rank condition)}} For every such choice,
		\[
		r(x_1 \vee \cdots \vee x_m) = \sum_{i=1}^{m} r(x_i).
		\]
	\end{itemize}
\end{defi}

\begin{example}\label{ex:poset-ideal-decomposition}
Figure~\ref{fig:poset-ideal-decomposition} gives a small ideal
decomposition \(\{I_1,I_2\}\). The two blue atoms generate \(I_1\), and
the two red atoms generate \(I_2\). Every nonzero element has a
nontrivial principal intersection with at least one of these ideals,
while every choice of one element from \(I_1\) and one from \(I_2\) has a
rank-additive join.
\begin{figure}[H]
\centering
\begin{tikzpicture}[
  scale=1,
  every node/.style={circle,draw,fill=black,inner sep=0.5pt,minimum size=1mm}]
  \node (bottom) at (-1,0) {};
  \node[fill=blue,minimum size=2mm] (b1) at (-4,1.5) {};
  \node[fill=blue,minimum size=2mm] (b2) at (-2,1.5) {};
  \node[fill=red,minimum size=2mm] (r1) at (0,1.5) {};
  \node[fill=red,minimum size=2mm] (r2) at (2,1.5) {};
  \node (u1) at (1,3) {};
  \node (u2) at (-1,3) {};
  \node (u3) at (-3,3) {};
  \draw[thick,blue] (bottom)--(b1);
  \draw[thick,blue] (bottom)--(b2);
  \draw[thick,red] (bottom)--(r1);
  \draw[thick,red] (bottom)--(r2);
  \draw (b1)--(u3);
  \draw (b1)--(u1);
  \draw (b2)--(u2);
  \draw (b2)--(u1);
  \draw (r1)--(u3);
  \draw (r1)--(u2);
  \draw (r2)--(u1);
\end{tikzpicture}
\caption{A meet-semilattice with the ideal decomposition
\(\{I_1,I_2\}\); the generators of \(I_1\) are blue and those of
\(I_2\) are red.}
\label{fig:poset-ideal-decomposition}
\end{figure}
\end{example}

Cartesian products give the basic model. If
\(L=L_1\times\cdots\times L_m\), the coordinate copies
\[
I_i=\{(\hat{0},\ldots,\hat{0},x_i,\hat{0},\ldots,\hat{0})
      :x_i\in L_i\}
\]
form an ideal decomposition of \(L\), although an ideal decomposition
need not arise from a poset product.

To prove Theorem \ref{thm:ideal-factorization}, we use the M\"{o}bius algebra introduced by
Solomon \cite{Solomon}; see also \cite{Greene,Stanley-book}. The
\emph{M\"{o}bius algebra} \(M(L)\) is the
\(\mathbb{Q}\)-vector space with basis \(L\) and multiplication
\[
x \cdot y =
\begin{cases}
x \vee y, & \text{if the join of $x$ and $y$ exists}; \\[3pt]
0, & \text{otherwise},
\end{cases}
\]
for all $x, y \in L$.
For each $x \in L$, define the element $\sigma_x \in M(L)$ by
\[
\sigma_x = \sum_{y \geq x}\mu(x, y)y.
\]
By the M\"{o}bius inversion formula \cite[Theorem 1.1, p.~399]{Stanley-book},
\begin{equation}\label{eq:mobius-inversion-basis}
	x = \sum_{y \geq x}\sigma_y,
\end{equation}
which implies that \(\{\sigma_x \mid x \in L\}\) forms a
\(\mathbb{Q}\)-basis for \(M(L)\).

\begin{lemma}\label{lem:mobius-factorization}
Let $L$ be a finite ranked meet-semilattice. If $L$ admits an ideal decomposition $\mathcal{I} = \{I_1, I_2, \ldots, I_m\}$,
then in the M\"{o}bius algebra $M(L)$,
	\[
		\sigma_{\hat{0}} = \prod_{i=1}^m \Big(\sum_{x_i \in I_i}\mu(x_i)x_i\Big).
	\]
	Equivalently, for each $x \in L$,
	\begin{equation}\label{eq:mobius-coefficient-factorization}
	\mu(x) =\sum_{\substack{(x_1, \ldots, x_m)\in I_1\times\cdots\times I_m\\[3pt]
x_{1}\vee \cdots \vee x_{m}=x}}\mu(x_{1})\cdots \mu(x_m).
	\end{equation}
\end{lemma}

\begin{proof}
By the join condition \hyperlink{join_condition}{$(\mathrm{I2})$},
\begin{align}\label{eq:mobius-product-expansion}
\prod_{i=1}^m \Big(\sum_{x_i \in I_i}\mu(x_i)x_i\Big)
=\sum_{(x_1, \ldots, x_m)\in I_1\times\cdots\times I_m}
\mu(x_1)\cdots\mu(x_m) x_1\vee\cdots \vee x_m.
\end{align}
Equation \eqref{eq:mobius-inversion-basis} gives
\[
x_1 \vee \cdots \vee x_m
= \sum_{y \geq x_1 \vee \cdots \vee x_m} \sigma_y.
\]
Substituting this into \eqref{eq:mobius-product-expansion},
we find
\begin{align*}
\prod_{i=1}^m \Big(\sum_{x_i \in I_i}\mu(x_i)x_i\Big)
=\sum_{(x_1,\ldots,x_m)\in I_1\times\cdots\times I_m} \mu(x_1)\cdots\mu(x_m) \sum_{y \geq x_1\vee \cdots \vee x_m}\sigma_y.
\end{align*}
Switching the order of summation and simplifying gives
\begin{align}\label{eq:mobius-sigma-expansion}
\prod_{i=1}^m \Big(\sum_{x_i \in I_i}\mu(x_i)x_i\Big)
=\sum_{y\in L}\sigma_y\prod_{i=1}^m\Big(\sum_{x_i\in I_i\cap \langle y \rangle }\mu(x_i)\Big).
\end{align}
For each \(y\in L\setminus\{\hat{0}\}\), the meet condition
\hyperlink{meet_condition}{$(\mathrm{I1})$} supplies an index \(j\)
such that
$I_j \cap \langle y \rangle=\langle z\rangle$ for some
\(z\neq\hat{0}\). Therefore
\[
\sum_{x_j\in I_j\cap \langle y \rangle }\mu(x_j)
=\sum_{x_j\leq z}\mu(\hat{0},x_j)=0.
\]
Thus, for any $y \in L \setminus \{\hat{0}\}$,
\[
\prod_{i=1}^m \Big(\sum_{x_i \in I_i \cap \langle y \rangle } \mu(x_i)\Big) = 0.
\]
Equation \eqref{eq:mobius-sigma-expansion} now reduces to
\begin{align}\label{eq:mobius-product-sigma0}
\prod_{i=1}^m \Big(\sum_{x_i \in I_i} \mu(x_i)x_i\Big) = \sigma_{\hat{0}}.
\end{align}
Using the definition of \(\sigma_{\hat{0}}\) and comparing the
coefficients of each basis element \(x\) on both sides of
\eqref{eq:mobius-product-sigma0} gives
\eqref{eq:mobius-coefficient-factorization}.
\end{proof}

\begin{proof}[Proof of Theorem \ref{thm:ideal-factorization}]
For each \(i\), choose \(z_i\in I_i\) with \(r(z_i)=r(I_i)\). Conditions
\hyperlink{join_condition}{$(\mathrm{I2})$} and
\hyperlink{rank_condition}{$(\mathrm{I3})$} give
\[
\sum_{i=1}^m r(I_i)=r(z_1\vee\cdots\vee z_m)\leq r(L),
\]
so \(c\geq0\).
\begin{align*}
	t^c\prod_{i=1}^m \chi(I_{i},t)
	&= t^c\prod_{i=1}^m
	\Big(\sum_{x_i\in I_i}\mu(x_i)t^{r(I_i)-r(x_i)}\Big) \\[6pt]
	&= \sum_{(x_1,\ldots,x_m)\in I_1\times\cdots\times I_m}
	\mu(x_1)\cdots\mu(x_m)
	t^{r(L)-r(x_1\vee\cdots\vee x_m)}
	\tag{by \hyperlink{rank_condition}{$(\mathrm{I3})$}} \\[6pt]
	&= \sum_{x\in L}
	\sum_{\substack{(x_1,\ldots,x_m)\in I_1\times\cdots\times I_m\\[3pt] x_{1}\vee \cdots \vee x_{m}=x}}
	\mu(x_{1})\cdots \mu(x_m) t^{r(L)-r(x)} \\[6pt]
	&= \sum_{x\in L} \mu(x)t^{r(L)-r(x)}
	\tag{by \eqref{eq:mobius-coefficient-factorization}} \\[6pt]
	&=\chi(L,t).
	\end{align*}
This completes the proof.
\end{proof}

\subsection{Intersection semilattices and nice partitions}
\label{subsec:arrangement-applications}

We now apply the preceding factorization to the intersection semilattice
\(L(\A)\). The nonvanishing of the M\"obius function on intersection
semilattices gives the following useful consequence.

\begin{lemma}\label{lem:join-generation}
If $L(\A)$ has an ideal decomposition
\(\mathcal{I}=\{I_1,I_2,\ldots,I_m\}\), then
\[
L(\A) = \{ X_1 \vee \cdots \vee X_m \mid X_i \in I_i, \;i = 1, \ldots, m \}.
\]
\end{lemma}

\begin{proof}
Suppose that some $X\in L(\A)$ cannot be expressed as such a join.
Then the index set in \eqref{eq:mobius-coefficient-factorization} is
empty, so Lemma \ref{lem:mobius-factorization} gives \(\mu(X)=0\). This
contradicts the sign theorem
\((-1)^{r(X)}\mu(X)>0\) for intersection semilattices; see
\cite[Theorem 3.10, p.~427]{Stanley-book}.
\end{proof}

For an essential arrangement \(\A\), \(r(L(\A))=d\) and
\(r(X)=d-\dim(X)\) for every \(X\in L(\A)\). Hence the arrangement
characteristic polynomial agrees with the poset characteristic polynomial
\eqref{eq:poset-characteristic-polynomial}:
\[
\chi_\A(t) = \sum_{X \in L(\A)} \mu(X) t^{\dim(X)}.
\]
Accordingly, we do not distinguish between these two expressions below.

\begin{coro}\label{cor:arrangement-characteristic-factorization}
If $L(\A)$ has an ideal decomposition $\mathcal{I} = \{I_1, I_2, \ldots, I_m\}$, then
	\[
	\chi_{\A}(t)= \prod_{i=1}^{m}\chi(I_{i},t).
	\]
\end{coro}
\begin{proof}
It suffices to show that \(c=0\) in Theorem
\ref{thm:ideal-factorization}. For each \(i\), choose
\(Y_i\in I_i\) with \(r(Y_i)=r(I_i)\). The join and rank conditions give
\[
\sum_{i=1}^{m} r(I_i) = \sum_{i=1}^{m} r(Y_i) = r\big(\bigvee_{i = 1}^m Y_i\big) \leq r(L(\A)) = d,
\]
where essentiality gives \(r(L(\A))=d\). For the reverse inequality, let
\(X\) be a maximal element of \(L(\A)\). By Lemma
\ref{lem:join-generation}, write
\(X=X_1\vee\cdots\vee X_m\) with \(X_i\in I_i\). Since every maximal
element has rank \(d\), the rank condition gives
\[
d = r(X) = \sum_{i=1}^{m} r(X_i) \leq \sum_{i=1}^{m} r(I_i).
\]
Consequently,
\begin{align}\label{eq:rank-sum-arrangement}
d = \sum_{i=1}^{m} r(I_i),
\end{align}
and hence \(c=0\).
\end{proof}

\begin{example}\label{ex:arrangement-ideal-decomposition}
Consider the arrangement \(\A\) in \(\mathbb{R}^3\) consisting of
\[
\begin{aligned}
H_1&\colon x-y=0,       & H_2&\colon x-y+1=0,
&H_3&\colon x=0,\\
H_4&\colon y=0,         & H_5&\colon z=0,
&H_6&\colon x-z=0.
\end{aligned}
\]
Let
\[
I_1=\langle H_1,H_2\rangle,
\qquad I_2=\langle H_3,H_4\rangle,
\qquad I_3=\langle H_5,H_6\rangle.
\]
The three ideals satisfy the meet, join, and rank conditions and hence
form an ideal decomposition of \(L(\A)\). Corollary
\ref{cor:arrangement-characteristic-factorization} yields
\[
\chi_{\A}(t)
=\chi(I_1,t)\chi(I_2,t)\chi(I_3,t)
=(t-2)^3.
\]
\end{example}

We next relate ideal decompositions to Terao's nice partitions
\cite{HT1992}. Let \(\pi\) be a partition of \(\A\), with blocks
\(\pi_1,\ldots,\pi_m\).
A \emph{section} of \(\pi\) is a tuple
\(S=(H_1,\ldots,H_p)\) whose hyperplanes lie in distinct blocks. It is
\emph{independent} if \(\bigcap_{i=1}^pH_i\neq\emptyset\) and
\(r(\bigcap_{i=1}^pH_i)=p\). The partition is \emph{independent} if every
section is independent.
For \(X\in L(\A)\), let
\[
\A_X:=\{H\in\A:X\subseteq H\}
\]
be the localization at \(X\). The partition \(\pi\) is \emph{nice} if it
is independent and, for every \(X\neq V\), some block \(\pi_i\) satisfies
\(|\pi_i\cap\A_X|=1\).

\begin{prop}\label{prop:nice-partition-ideal-decomposition}
Let \(\pi=(\pi_1,\ldots,\pi_m)\) be a partition of \(\A\). Then \(\pi\)
is nice if and only if
\(\{\langle\pi_1\rangle,\ldots,\langle\pi_m\rangle\}\) is an ideal
decomposition of \(L(\A)\).
\end{prop}

\begin{proof}
Put $I_i=\langle \pi_i\rangle$ for $1\leq i\leq m$. Since each element of
\(\pi_i\) is a hyperplane, \(I_i=\{V\}\cup \pi_i\).

Assume first that \(\pi\) is nice. Let \(X\in L(\A)\setminus\{V\}\). By the
singleton condition, there is an index \(i\) such that
\(\pi_i\cap \A_X=\{H\}\). Hence
\[
I_i\cap \langle X\rangle=\{V,H\}=\langle H\rangle,
\]
which is a nontrivial principal ideal. This proves the meet condition.
For the join and rank conditions, choose \(X_i\in I_i\). After deleting the
entries equal to \(V\), the remaining elements form a section of \(\pi\).
The independence of \(\pi\) says exactly that their intersection is
nonempty and has rank equal to the number of chosen hyperplanes. Thus
\(X_1\vee\cdots\vee X_m\) exists and its rank is
\(\sum_i r(X_i)\). Therefore \(\{I_1,\ldots,I_m\}\) is an ideal
decomposition.

For the converse, suppose that \(\{I_1,\ldots,I_m\}\) is an ideal
decomposition. Let \(S\) be a section of \(\pi\), say
\[
S=(H_1,\ldots,H_p),\qquad
H_j\in\pi_{i_j},\quad i_1<\cdots<i_p.
\]
Apply the join and rank conditions to these \(p\) hyperplanes, and choose
\(V\) in all remaining positions. Then \(H_1\cap\cdots\cap H_p\) is
nonempty and has rank \(p\).
Hence \(\pi\) is independent. Finally, for
\(X\in L(\A)\setminus\{V\}\), the meet condition gives an index \(i\) such
that \(I_i\cap \langle X\rangle\) is a nontrivial principal ideal. Since
\(I_i=\{V\}\cup\pi_i\), this is possible only when
\(\pi_i\cap \A_X\) consists of exactly one hyperplane. Thus \(\pi\) is
nice.
\end{proof}

For each \(i\in\{1,\ldots,m\}\),
\(\chi(\langle\pi_i\rangle,t)=t-|\pi_i|\).
Combining this with Corollary
\ref{cor:arrangement-characteristic-factorization} and Proposition
\ref{prop:nice-partition-ideal-decomposition}, we recover the following
corollary, originally obtained from the graded factorization of the
Orlik--Solomon algebra \cite[Corollary 2.9]{HT1992}.
\begin{coro}\label{cor:nice-characteristic-factorization}
	If $\A$ admits a nice partition $\pi = (\pi_{1},\ldots,\pi_{m})$, then
	\[
	\chi_{\A}(t) = \prod_{i=1}^{m}(t - |\pi_{i}|).
	\]
\end{coro}

\subsection{Semimatroids and geometric semilattices}

The flat semilattices of semimatroids are precisely the geometric
semilattices \cite{Ardila-semimatroids}. Theorem
\ref{thm:ideal-factorization} is therefore a result about ranked
meet-semilattices themselves and applies directly to geometric
semilattices, not only to those arising from arrangements.
Every lower interval of a geometric semilattice is a geometric lattice,
so its M\"obius values have alternating signs. The argument of Lemma
\ref{lem:join-generation} then shows that the rank correction vanishes.
Hence every ideal decomposition
\(\{I_1,\ldots,I_m\}\) of \(L(\mathfrak S)\) satisfies
\[
\chi(L(\mathfrak S),t)=\prod_{i=1}^m\chi(I_i,t).
\]

In the join-closed case, ideal decompositions of simple semimatroids
in fact come from direct sums. Since the atoms of \(L(\mathfrak S)\)
correspond to the ground-set elements, the ideals determine the subsets
\(E_i\) in the following theorem.

\begin{theorem}\label{thm:join-closed-semimatroid-direct-sum}
Let \(\mathfrak S=(E,\mathcal C,r)\) be a finite simple semimatroid, and
suppose that \(L(\mathfrak S)\) has an ideal decomposition
\(\mathcal I=\{I_1,\ldots,I_m\}\). Assume that each \(I_i\) is closed
under every join that exists in \(L(\mathfrak S)\). Define
\[
E_i:=\{e\in E:\operatorname{cl}(\{e\})\in I_i\},
\qquad
\mathfrak S_i:=\mathfrak S|_{E_i}.
\]
Then
\[
E=E_1\sqcup\cdots\sqcup E_m,
\qquad
L(\mathfrak S_i)\cong I_i,
\]
and
\[
\mathfrak S\cong
\mathfrak S_1\oplus\cdots\oplus\mathfrak S_m.
\]
In particular, the join map
\[
I_1\times\cdots\times I_m\longrightarrow L(\mathfrak S)
\]
is a rank-preserving poset isomorphism.
\end{theorem}

\begin{proof}
For \(e\in E\), let \(a_e=\operatorname{cl}(\{e\})\), an atom of
\(L(\mathfrak S)\). Applying the meet condition
\hyperlink{meet_condition}{$(\mathrm{I1})$} to \(a_e\) shows that
\(a_e\in I_i\) for some \(i\). It cannot belong to two distinct ideals:
if \(a_e\in I_i\cap I_j\) with \(i\neq j\), choose \(a_e\) in the
\(i\)-th and \(j\)-th positions and \(\hat{0}\) elsewhere. The rank
condition would give
\[
1=r(a_e\vee a_e)=r(a_e)+r(a_e)=2,
\]
a contradiction. Hence
\(E=E_1\sqcup\cdots\sqcup E_m\).

Every flat in \(I_i\) is the join of the atoms below it, and all those
atoms are indexed by \(E_i\). Conversely, if \(A_i\subseteq E_i\) is
central, then \(\operatorname{cl}_{\mathfrak S}(A_i)\) is the existing
join of the atoms indexed by \(A_i\). Join-closedness therefore places
this flat in \(I_i\). No atom indexed by \(E\setminus E_i\) can lie below
this closure: otherwise the ideal \(I_i\) would contain an atom assigned
to a different block. Thus closure in \(\mathfrak S\) introduces no
element outside \(E_i\) and agrees with closure in the restriction
\(\mathfrak S_i\). These observations identify
\(L(\mathfrak S_i)\) with \(I_i\).

To identify the semimatroid, take subsets \(A_i\subseteq E_i\). If each
\(A_i\) is central in \(\mathfrak S_i\), then
\(X_i=\operatorname{cl}(A_i)\in I_i\), and the join condition gives
\(X_1\vee\cdots\vee X_m\). In a flat semilattice, this is equivalent to
the existence of a central set containing \(A_1\cup\cdots\cup A_m\);
since central sets form a simplicial complex, the union itself is
central. Conversely, if \(A\subseteq E\) is central, then each
\(A_i=A\cap E_i\) is central and the join of the flats
\(\operatorname{cl}(A_i)\) is \(\operatorname{cl}(A)\). In both
directions, the rank condition gives
\[
r(A_1\cup\cdots\cup A_m)=\sum_{i=1}^m r(A_i),
\qquad
r(A)=\sum_{i=1}^m r(A\cap E_i).
\]
Thus the central sets and ranks agree with those of
\(\mathfrak S_1\oplus\cdots\oplus\mathfrak S_m\).

Finally, the flat semilattice of a direct sum is the direct product of
the flat semilattices of its summands. Under the identifications above,
this product isomorphism is exactly the join map.
\end{proof}

\section{Modular ideals and geometric realization}
\label{sec:modular-realization}

\subsection{Modular ideals and two-factor decompositions}
Stanley's modular-element factorization for geometric lattices
\cite[Theorem 4.13, p.~439]{Stanley-book} applies directly to central
arrangements. For affine arrangements, \(L(\A)\) need only be a
meet-semilattice, so we use the following ideal-theoretic replacement.

\begin{defi}\label{def:modular-ideal}
Let \(\A\) be an affine arrangement and \(I\) an ideal of \(L(\A)\).
Define
\begin{align}\label{eq:orthogonal-ideal}
I^\perp = \{Y \in L(\A) \mid Y \wedge X = \hat{0}
\text{ for all } X \in I\}.
\end{align}
This set is an ideal, called the \emph{orthogonal ideal} of \(I\). The
ideal \(I\) is \emph{modular} if
\(\{I,I^\perp\}\) forms an ideal decomposition of \(L(\A)\).
\end{defi}

\begin{example}\label{ex:modular-ideal-realization}
Let \(\A\) be the affine arrangement in \(\mathbb{R}^3\) consisting of
the following seven hyperplanes:
\begin{align*}
    H_1 &\colon x + y - 1 = 0, \quad &
    H_2 &\colon x = 0, \quad &
    H_3 &\colon x - 1 = 0, \quad &
    H_4 &\colon x - y = 0, \\[6pt]
    H_5 &\colon y + z - 1 = 0, \quad &
    H_6 &\colon y = 0, \quad &
    H_7 &\colon y - 1 = 0.
\end{align*}
Consider two ideals of \(L(\A)\):
\begin{align*}
I_1 & = \langle H_1, H_2, H_3 \rangle \\[6pt]
I_2 & = \langle H_4 \cap H_5,\ H_5 \cap H_6,\ H_5 \cap H_7 \rangle.
\end{align*}
Write \(L_{45}=H_4\cap H_5\), \(L_{56}=H_5\cap H_6\), and
\(L_{57}=H_5\cap H_7\). The maximal elements of \(I_1\) are
\(H_1,H_2,H_3\), while those of \(I_2\) are
\(L_{45},L_{56},L_{57}\). The following table lists the points
\(H_i\cap L_{ab}\), where \(1\leq i\leq3\) and
\(L_{ab}\in\{L_{45},L_{56},L_{57}\}\):
\[
\begin{array}{c|ccc}
 &L_{45}&L_{56}&L_{57}\\ \hline
H_1&(1/2,1/2,1/2)&(1,0,1)&(0,1,0)\\
H_2&(0,0,1)&(0,0,1)&(0,1,0)\\
H_3&(1,1,0)&(1,0,1)&(1,1,0).
\end{array}
\]
The table verifies the case \(X_1=H_i\) and \(X_2=L_{ab}\). If
\(X_1=\hat{0}\) or \(X_2=\hat{0}\), there is nothing to check. Otherwise \(X_2\)
is either one of the lines \(L_{ab}\), already
covered by the table, or a hyperplane containing one of them. In the latter
case \(H_i\cap X_2\) contains the corresponding table point and has rank
\(2=1+1\). Thus every join \(X_1\vee X_2\), with \(X_i\in I_i\), exists and
satisfies \(r(X_1\vee X_2)=r(X_1)+r(X_2)\).

Figure~\ref{fig:modular-ideals-example} shows the full intersection
semilattice. Blue and red nodes denote the nonzero elements of \(I_1\) and
\(I_2\), respectively; gray nodes belong to neither ideal. We now verify the
meet condition. If \(\langle X\rangle\cap I_1\) is
nontrivial principal, there is nothing to prove. From
Figure~\ref{fig:modular-ideals-example}, the remaining possibilities are
exactly the following: either \(X\in I_2\), or the localization \(\A_X\) is
one of
\[
\{H_1,H_2,H_7\},\quad
\{H_1,H_3,H_6\},\quad
\{H_1,H_2,H_5,H_7\},\quad
\{H_1,H_3,H_5,H_6\}.
\]
If \(X\in I_2\), then \(\langle X\rangle\cap I_2=\langle X\rangle\). In the
four listed cases, \(\langle X\rangle\cap I_2\) is generated, respectively,
by \(H_7,H_6,L_{57},L_{56}\). Hence the meet condition holds, and
\(\{I_1,I_2\}\) satisfies the three axioms of Definition
\ref{def:ideal-decomposition}.
Thus \(I_1\) and \(I_2\) are mutually orthogonal modular ideals.

\begin{figure}[H]
\centering
\resizebox{0.98\textwidth}{!}{%
\begin{tikzpicture}[
  x=1cm,y=1cm,
  posetnode/.style={draw,rounded corners=1.2pt,inner sep=1.5pt,
    minimum height=4.8mm,font=\scriptsize},
  bluev/.style={posetnode,fill=blue!16,draw=blue!70!black},
  redv/.style={posetnode,fill=red!16,draw=red!70!black},
  grayv/.style={posetnode,fill=gray!12,draw=gray!65!black},
  bottomv/.style={posetnode,fill=gray!8,draw=gray!65!black},
  cover/.style={draw=gray!58,line width=0.36pt}]
  \node[bottomv] (V) at (0,0) {$V$};
  \node[bluev] (H1) at (-5.4,1.35) {$H_1$};
  \node[bluev] (H2) at (-3.6,1.35) {$H_2$};
  \node[bluev] (H3) at (-1.8,1.35) {$H_3$};
  \node[redv]  (H4) at (0,1.35) {$H_4$};
  \node[redv]  (H5) at (1.8,1.35) {$H_5$};
  \node[redv]  (H6) at (3.6,1.35) {$H_6$};
  \node[redv]  (H7) at (5.4,1.35) {$H_7$};
  \node[grayv] (R127) at (-6.5,3.05) {$127$};
  \node[grayv] (R136) at (-5.2,3.05) {$136$};
  \node[grayv] (R14)  at (-3.9,3.05) {$14$};
  \node[grayv] (R15)  at (-2.6,3.05) {$15$};
  \node[grayv] (R25)  at (-1.3,3.05) {$25$};
  \node[grayv] (R35)  at (0,3.05) {$35$};
  \node[redv]  (R45)  at (1.3,3.05) {$45$};
  \node[redv]  (R56)  at (2.6,3.05) {$56$};
  \node[redv]  (R57)  at (3.9,3.05) {$57$};
  \node[grayv] (R246) at (5.2,3.05) {$246$};
  \node[grayv] (R347) at (6.5,3.05) {$347$};
  \node[grayv] (P1257) at (-4.8,4.75) {$1257$};
  \node[grayv] (P1356) at (-2.4,4.75) {$1356$};
  \node[grayv] (P145)  at (0,4.75) {$145$};
  \node[grayv] (P2456) at (2.4,4.75) {$2456$};
  \node[grayv] (P3457) at (4.8,4.75) {$3457$};
  \foreach \u/\v in {V/H1,V/H2,V/H3,V/H4,V/H5,V/H6,V/H7,
    H1/R127,H1/R136,H1/R14,H1/R15,
    H2/R127,H2/R246,H2/R25,
    H3/R136,H3/R347,H3/R35,
    H4/R14,H4/R246,H4/R347,H4/R45,
    H5/R15,H5/R25,H5/R35,H5/R45,H5/R56,H5/R57,
    H6/R136,H6/R246,H6/R56,
    H7/R127,H7/R347,H7/R57,
    R127/P1257,R136/P1356,R14/P145,
    R15/P1257,R15/P1356,R15/P145,
    R25/P1257,R25/P2456,R35/P1356,R35/P3457,
    R45/P145,R45/P2456,R45/P3457,
    R56/P1356,R56/P2456,R57/P1257,R57/P3457,
    R246/P2456,R347/P3457}
    {\draw[cover] (\u)--(\v);}
\end{tikzpicture}}
\caption{The intersection semilattice in Example
\ref{ex:modular-ideal-realization}. A label such as \(127\) denotes
\(H_1\cap H_2\cap H_7\).}
\label{fig:modular-ideals-example}
\end{figure}

The subarrangements \(\{H_1,H_2,H_3\}\) and
\(\{H_4,H_5,H_6,H_7\}\) do not realize the two ideals: for example,
\(H_1\cap H_2,H_1\cap H_3\in L(\{H_1,H_2,H_3\})\setminus I_1\).
The translated restrictions constructed in Theorem
\ref{thm:geometric-realization} instead have the combinatorial forms
shown in Figure~\ref{fig:geometric-realizations-example}.
\begin{figure}[H]
\centering
\begin{tikzpicture}[scale=1.05,every node/.style={circle}]
  \draw[thick] (4.7,1)--(8.4,1);
  \draw[thick] (4.7,0)--(8.4,0);
  \draw[thick] (4.7,-1)--(8.4,-1);
  \draw[thick] (5.1,1.45)--(8.1,-1.45);
  \node[fill=red,circle,inner sep=2pt] at (5.55,1) {};
  \node[fill=red,circle,inner sep=2pt] at (6.6,0) {};
  \node[fill=red,circle,inner sep=2pt] at (7.65,-1) {};
  \draw[thick] (-2,0)--(2,0);
  \node[fill=blue,circle,inner sep=2pt] at (-1.5,0) {};
  \node[fill=blue,circle,inner sep=2pt] at (0,0) {};
  \node[fill=blue,circle,inner sep=2pt] at (1.5,0) {};
\end{tikzpicture}
\caption{Geometric realizations of \(I_1\) (left) and \(I_2\) (right)
from Example~\ref{ex:modular-ideal-realization}.}
\label{fig:geometric-realizations-example}
\end{figure}

The example also shows that modular ideals need not be join-closed. Indeed,
\(H_1,H_2\in I_1\), but their join exists and satisfies
\[
H_1\vee H_2=H_1\cap H_2=H_1\cap H_2\cap H_7\notin I_1.
\]
\end{example}

For a central arrangement \(\A\), an element \(Z\in L(\A)\) is
\emph{modular} if
\[
r(Z)+r(X)=r(Z\wedge X)+r(Z\vee X)
\]
for every \(X\in L(\A)\).
For central arrangements, modular ideals recover the usual modular elements,
up to taking orthogonal ideals.

\begin{prop}\label{prop:central-modular-ideal}
Let $\A$ be a central arrangement. We have
\begin{enumerate}
    \item[{\rm (1)}] If \(Z \in L(\A)\) is a modular element, then \(\langle Z \rangle\) is a modular ideal.
    \item[{\rm (2)}] If $I$ is a modular ideal of $L(\A)$, then either $I$ or $I^\perp$ is a principal ideal generated by a modular element.
\end{enumerate}
\end{prop}

\begin{proof}
(1) Assume that \(Z\) is modular in \(L(\A)\), so
\begin{align}\label{eq:modular-element-rank}
r(Z) + r(X) = r(Z \wedge X) + r(Z \vee X)
\end{align}
for every \(X\in L(\A)\). We verify that
\(\{\langle Z\rangle,\langle Z\rangle^\perp\}\) satisfies the
three conditions in Definition \ref{def:ideal-decomposition}.
\begin{enumerate}
    \item[\hyperlink{meet_condition}{$(\mathrm{I1})$}]
    Fix \(X \in L(\A) \setminus \{\hat{0}\}\).
    If $X \wedge Z \neq \hat{0}$, then
    $\langle X \wedge Z \rangle = \langle Z \rangle \cap \langle X \rangle$
    is the nontrivial principal ideal we need. If \(X\wedge Z=\hat{0}\),
    then for every \(Z'\leq Z\) we also have \(X\wedge Z'=\hat{0}\);
    indeed, \(X\wedge Z'\leq X\wedge Z=\hat{0}\). Thus
    \(X\in\langle Z\rangle^\perp\), and
    \(\langle X\rangle=\langle X\rangle\cap
    \langle Z\rangle^\perp\) is the required nontrivial principal
    ideal.
    \item[\hyperlink{join_condition}{$(\mathrm{I2})$}] Since $\A$ is central, $L(\A)$ is a lattice, so the join condition holds.
    \item[\hyperlink{rank_condition}{$(\mathrm{I3})$}]
For any \(X\in\langle Z\rangle^\perp\),
\(Z\wedge X=V=\hat{0}\), and hence
\(r(Z)+r(X)=r(Z\vee X)\) by
\eqref{eq:modular-element-rank}. In terms of dimensions, this equality is
\[
\dim(Z) + \dim(X) = d + \dim(Z \cap X),
\]
so \(Z+X=V\). If \(Z'\in\langle Z\rangle\), then \(Z\subseteq Z'\)
geometrically, and therefore \(Z'+X=V\) as well. Applying the dimension
formula to \(Z'\) and \(X\) gives
\(r(Z'\vee X)=r(Z')+r(X)\), which is the required rank condition.
\end{enumerate}
Thus \(\langle Z\rangle\) is a modular ideal.

(2) Assume that $I$ is a modular ideal of $L(\A)$, that is, $\{I, I^\perp\}$ is an ideal decomposition of $L(\A)$.
Applying the meet condition to \(\hat{1}\), we find that either \(I\) or
\(I^\perp\) is a nontrivial principal ideal. Interchanging the
two ideals if necessary, write \(I=\langle Z\rangle\).
It remains to prove that \(Z\) is modular.
From the rank condition \hyperlink{rank_condition}{$(\mathrm{I3})$},
it follows that for any \(X\in\langle Z\rangle^\perp\), equivalently
for any \(X\) with \(X\wedge Z=\hat{0}\),
\(r(X)+r(Z)=r(X\vee Z)\). Brylawski's rank criterion
\cite[Theorem~3.3]{Brylawski} now implies that \(Z\) is modular.
\end{proof}

We next compare the modular-ideal factorization with Stanley's classical
factorization of the characteristic polynomial. Let \(\A\) be a central
arrangement. For \(Z\in L(\A)\), write
\[
\A_Z:=\{H\in\A\mid Z\subseteq H\}
\]
for the localization at \(Z\).
When \(\A_Z\) is viewed as an arrangement in the original ambient space \(V\),
its characteristic polynomial is
\[
\chi_{\A_Z}(t)=t^{d-r(Z)}\chi(\langle Z\rangle,t).
\]
Here \(\chi(\langle Z\rangle,t)\) denotes the poset characteristic polynomial
of the principal ideal \(\langle Z\rangle\).
Stanley \cite{Stanley1971} proved that if $Z\in L(\A)$ is a modular element, then
\begin{align}\label{eq:stanley-factorization}
\chi_\A(t) = \chi(\langle Z\rangle,t)
\Big( \sum_{Y \wedge Z = \hat{0}} \mu(Y)t^{d - r(Y) - r(Z)} \Big).
\end{align}
For a modular ideal \(I\subseteq L(\A)\), Corollary
\ref{cor:arrangement-characteristic-factorization} gives
\begin{align}\label{eq:central-characteristic-factorization}
\chi_\A(t) = \chi(I, t) \chi(I^\perp, t).
\end{align}
When \(Z\) is modular and \(I=\langle Z\rangle\), Proposition
\ref{prop:central-modular-ideal} and the definition of \(I^\perp\) identify
this formula with \eqref{eq:stanley-factorization}. Over an infinite field,
the next subsection gives geometric realizations of these factors.

\subsection{Geometric realization of modular ideals}
\label{subsec:geometric-realization}
Choose a maximal \(Y\in I^\perp\), translate it to avoid finitely many
proper affine subspaces, and restrict \(\A_I\) to that translate. We prove
that the resulting intersection semilattice is \(I\).

Let \(\A\) be an affine arrangement. For every subarrangement
\(\B\subseteq\A\), regard \(L(\B)\) as a subposet of \(L(\A)\). For an
ideal \(I\subseteq L(\A)\), an affine subspace \(Y\subseteq V\), and
\(\alpha\in V\), define
\[
(\A_I)^{Y+\alpha}
:=\{H\cap(Y+\alpha):H\in\A_I,\ 
H\cap(Y+\alpha)\notin\{\emptyset,Y+\alpha\}\},
\]
an arrangement in \(Y+\alpha\). We seek \(Y\) and \(\alpha\) such that
\[
I \cong L\left((\A_I)^{Y + \alpha}\right)
\]
for the given modular ideal \(I\). We use the following elementary
properties.

\begin{prop}\label{prop:orthogonal-ideal-symmetry}
	Let \(I_1\) and \(I_2\) be ideals of \(L(\A)\). Then the following statements are equivalent:
	\begin{enumerate}
		\item[{\rm (1)}] $\{I_1, I_2\}$ forms an ideal decomposition of $L(\A)$;
		\item[{\rm (2)}] $I_1$ is a modular ideal of $L(\A)$ and $I_1^\perp = I_2$;
		\item[{\rm (3)}] $I_2$ is a modular ideal of $L(\A)$ and $I_2^\perp = I_1$.
	\end{enumerate}
\end{prop}
\begin{proof}
By symmetry, it suffices to prove that (1) is equivalent to (2).
If \(I_1\) is modular and \(I_1^\perp=I_2\), then \((1)\) follows
from Definition \ref{def:modular-ideal}. Conversely, if \(\{I_1,I_2\}\) is an ideal
decomposition of \(L(\A)\), it remains to prove \(I_1^\perp=I_2\).
	Let \(Y\in I_2\) and \(X\in I_1\). Since \(I_1\) is an ideal,
	\(X\wedge Y\in I_1\). Apply the rank condition to the pair
	\(X\wedge Y\in I_1\) and \(Y\in I_2\). Their join is \(Y\), so
	\[
	r(X \wedge Y) + r(Y) = r(Y),
	\]
	which implies \(X\wedge Y=\hat{0}\). Since this holds for every
	\(X\in I_1\), we have \(Y\in I_1^\perp\).
	On the other hand, for any \(Y\in I_1^\perp\), Lemma~\ref{lem:join-generation}
	guarantees the existence of \(Y_1\in I_1\) and \(Y_2\in I_2\) such that
	\[
	Y = Y_1 \vee Y_2.
	\]
	Since \(Y\in I_1^\perp\) and \(Y_1\leq Y\), we have
	\(Y_1=Y_1\wedge Y=\hat{0}\). Hence \(Y=Y_2\in I_2\).
	Thus \(I_2 = I_1^\perp\).
\end{proof}

\begin{coro}\label{cor:orthogonal-involution}
If \(I\) is a modular ideal of \(L(\A)\), then \(I^\perp\) is also
modular and \((I^\perp)^\perp=I\).
\end{coro}

\begin{remark*}
The modularity hypothesis in Corollary
\ref{cor:orthogonal-involution} is essential. Even when
\(I^\perp\) is modular, the equality
\((I^\perp)^\perp=I\) may fail if \(I\) itself is not
modular. For example, consider the braid arrangement
\[
\B_4 = \{H_{ij}: x_i - x_j = 0 \mid 1 \leq i < j \leq 4\}
\]
in $\mathbb{R}^4$. Take the non-modular ideal
$
I = \{\mathbb{R}^4, H_{12}, H_{13}, H_{23}\}.
$
A direct calculation gives
\begin{align*}
I^\perp&=\{\mathbb{R}^{4}, H_{14}, H_{24}, H_{34}\},\\[6pt]
(I^\perp)^\perp&=\{\mathbb{R}^{4}, H_{12}, H_{13}, H_{23}, H_{12}\cap H_{13}\cap H_{23}\}.
\end{align*}
Indeed, put
\(Z=H_{12}\cap H_{13}\cap H_{23}\) and
\(J=\langle Z\rangle=(I^\perp)^\perp\). In the partition lattice
\(L(\B_4)\cong\Pi_4\), the pair \(\{I^\perp,J\}\) is an ideal
decomposition. Hence \(I^\perp\) is
modular, whereas \(I\neq(I^\perp)^\perp\).
\end{remark*}	

\begin{prop}\label{prop:graded-modular-ideals}
All modular ideals of $L(\A)$ are graded.
\end{prop}

\begin{proof}
Let \(I\) be a modular ideal of \(L(\A)\), and let \(X_1\) be an
arbitrary maximal element of \(I\). We prove that \(r(X_1)=r(I)\).
Choose \(X_2\in I^\perp\) with
\(r(X_2)=r(I^\perp)\).
The join condition guarantees that \(X:=X_1\vee X_2\) exists. The rank
condition and \eqref{eq:rank-sum-arrangement} give
\[
r(X) = r(X_1) + r(X_2) \leq r(I) + r(I^\perp) = r(L(\A)).
\]
Extend \(X\) to a maximal element \(Y\) of \(L(\A)\). Since \(L(\A)\)
is graded, \(r(Y)=r(L(\A))\). By Lemma
\ref{lem:join-generation}, write \(Y=Y_1\vee Y_2\), where
\(Y_1\in I\) and \(Y_2\in I^\perp\). Then
\[
r(L(\A))=r(Y) = r(Y_1) + r(Y_2) \leq r(I) + r(I^\perp) = r(L(\A)),
\]
so \(r(Y_1)=r(I)\) and
\(r(Y_2)=r(I^\perp)\).
In particular, \(Y_1\) is maximal in \(I\) and \(Y_2\) is maximal in
\(I^\perp\), since rank strictly increases between distinct
comparable elements. The meet condition says that one of
\(\langle Y\rangle\cap I\) and
\(\langle Y\rangle\cap I^\perp\) is nontrivial principal.
\begin{enumerate}
\item[(1)] Suppose that
\(\langle Y\rangle\cap I=\langle Z\rangle\), where
\(Z\neq\hat{0}\). Since
\[
X_1\leq X=X_1\vee X_2\leq Y,
\]
we have
\(X_1\in\langle Y\rangle\cap I=\langle Z\rangle\), and hence
\(X_1\leq Z\). Likewise,
\(Y_1\leq Y\) and \(Y_1\in I\) imply
\(Y_1\in\langle Y\rangle\cap I=\langle Z\rangle\), so
\(Y_1\leq Z\). Since \(Z\in I\), the maximality of \(X_1\) and
\(Y_1\) in \(I\) therefore gives \(X_1=Y_1=Z\), and hence
\(r(X_1)=r(I)\).
\item[(2)] Suppose that
\(\langle Y\rangle\cap I^\perp\) is nontrivial principal.
Applying the argument from the first case to the maximal elements \(X_2,Y_2\)
of \(I^\perp\) gives
\[
X_2=Y_2,\qquad
\langle Y\rangle\cap I^\perp=\langle Y_2\rangle.
\]
We claim that \(X_1\vee Y_2=Y\). Once this is proved, the rank condition
applied to \(X_1\in I\) and \(Y_2\in I^\perp\) gives
$
r(Y) = r(X_1) + r(Y_2),
$
and comparison with \(r(Y)=r(Y_1)+r(Y_2)\) yields
\(r(X_1)=r(Y_1)=r(I)\).

To prove the claim, we first establish
\begin{align}\label{eq:graded-modular-case-ii}
\langle Y \rangle \cap I = \langle Y \rangle \cap \langle Y_2 \rangle^\perp.
\end{align}
Since \(\langle Y_2\rangle\subseteq I^\perp\), Corollary
\ref{cor:orthogonal-involution} gives
\[
I=(I^\perp)^\perp
\subseteq\langle Y_2\rangle^\perp.
\]
This proves one inclusion in \eqref{eq:graded-modular-case-ii}. Conversely,
take \(Z\in\langle Y\rangle\cap\langle Y_2\rangle^\perp\).
By Lemma \ref{lem:join-generation},
we can write \(Z=Z_1\vee Z_2\) for some \(Z_1\in I\) and
\(Z_2\in I^\perp\).
Since \(Z_2\leq Z\leq Y\),
\[
Z_2\in\langle Y\rangle\cap I^\perp
=\langle Y_2\rangle.
\]
Thus \(Z_2\leq Y_2\). Together with \(Z_2\leq Z\) and
\(Z\wedge Y_2=\hat{0}\), this gives \(Z_2=\hat{0}\), and hence
\(Z=Z_1\in I\cap\langle Y\rangle\). This proves
\eqref{eq:graded-modular-case-ii}.

We next show that \(Y_2\) is modular in
\(L(\A_Y)=\langle Y\rangle=[\hat{0},Y]\). This interval is a
geometric lattice with the rank function inherited from \(L(\A)\).
For elements of the interval, its meet agrees with the meet in
\(L(\A)\), and their join exists in the interval because \(Y\) is a
common upper bound; this join also agrees with their join in
\(L(\A)\). Consequently,
\[
\{Z\in\langle Y\rangle:Z\wedge Y_2=\hat{0}\}
=\langle Y\rangle\cap\langle Y_2\rangle^\perp.
\]
If
\(Z\in\langle Y\rangle\cap\langle Y_2\rangle^\perp\), then
\eqref{eq:graded-modular-case-ii} gives \(Z\in I\). The rank condition
therefore yields
$
r(Z) + r(Y_2) = r(Z \vee Y_2).
$
Thus the required rank equality holds for every element of
\(\langle Y\rangle\) whose meet with \(Y_2\) is \(\hat{0}\).
Brylawski's rank criterion \cite[Theorem~3.3]{Brylawski} therefore shows that
\(Y_2\) is modular in \(\langle Y\rangle\).
Finally, suppose for contradiction that \(X_1\vee Y_2<Y\).
The interval \(\langle Y\rangle\) is a geometric lattice and is
therefore atomic. Hence there is a hyperplane
\(H\in\langle Y\rangle\) such that
\(H\not\leq X_1\vee Y_2\); otherwise \(X_1\vee Y_2\) would be the join
of all atoms below \(Y\), and hence would equal \(Y\).
Consequently,
\begin{align*}
r(H \vee X_1 \vee Y_2)
&= 1 + r(X_1 \vee Y_2) \tag{by $H \not\leq X_1 \vee Y_2$}\\[6pt]
&= 1 + r(X_1) + r(Y_2)\tag{by the rank condition \hyperlink{rank_condition}{$(\mathrm{I3})$}} \\[6pt]
&= r(H \vee X_1) + r(Y_2). \tag{by $H \not\leq X_1$}
\end{align*}
On the other hand, the modularity of $Y_2$ gives
\[
r(H \vee X_1 \vee Y_2) = r(H \vee X_1) + r(Y_2) - r((H \vee X_1) \wedge Y_2).
\]
Combining the two identities above,
we conclude that $(H \vee X_1) \wedge Y_2 = \hat{0}$, which implies
$H \vee X_1 \in \langle Y \rangle \cap \langle Y_2 \rangle^\perp$.
By \eqref{eq:graded-modular-case-ii}, $H \vee X_1 \in \langle Y \rangle \cap I$,
which contradicts the maximality of $X_1$.
Thus \(X_1\vee Y_2=Y\).
\end{enumerate}
Since \(X_1\) was arbitrary, every maximal element of \(I\) has rank
\(r(I)\). Thus \(I\) is graded.
\end{proof}

\begin{lemma}\label{lem:avoidance-for-realization}
Let \(I\) be a modular ideal of \(L(\A)\), and let \(Y\) be a maximal
element of \(I^\perp\). If \(X\in L(\A_I)\setminus I\), then
either \(X\vee Y\) does not exist or \(X\wedge Y\neq\hat{0}\).
\end{lemma}
\begin{proof}
Fix \(X\in L(\A_I)\setminus I\). We first claim that, for every
\(X'\geq X\), the ideal \(\langle X'\rangle\cap I\) is not nontrivial
principal. Suppose instead that
\(\langle X'\rangle\cap I=\langle Z\rangle\) for some
\(Z\neq\hat{0}\). Since \(X\in L(\A_I)\), it is the join of all
hyperplanes in \(\langle X\rangle\cap\A_I\). Moreover,
\[
\langle X \rangle \cap \A_I = \langle X \rangle \cap I \cap \A \subseteq \langle X' \rangle \cap I = \langle Z \rangle.
\]
Taking their join gives \(X\leq Z\). On the other hand,
\(\langle X'\rangle\cap I=\langle Z\rangle\) and \(X\leq X'\) give
\[
\langle X \rangle \cap I
=\langle X\rangle\cap\langle Z\rangle
=\langle X\wedge Z\rangle.
\]
Since \(X\notin I\), the left-hand side does not contain \(X\), so
\(X\neq X\wedge Z\). This contradicts \(X\leq Z\), which implies
\(X=X\wedge Z\), and proves the claim.

Suppose, to the contrary, that \(X\vee Y\) exists and
\(X\wedge Y=\hat{0}\). By the claim,
\(\langle X\vee Y\rangle\cap I\) is not nontrivial principal. The meet
condition therefore gives
\[
\langle X\vee Y\rangle\cap I^\perp=\langle Z'\rangle
\]
for some \(Z'\in I^\perp\setminus\{\hat{0}\}\). Since \(Y\)
belongs to the left-hand side, \(Y\leq Z'\). The maximality of \(Y\) in
\(I^\perp\) yields \(Y=Z'\), and hence
\[
\langle X\vee Y\rangle\cap I^\perp=\langle Y\rangle.
\]
Lemma \ref{lem:join-generation} ensures that
\(X=X_1\vee X_2\) for some \(X_1\in I\) and
\(X_2\in I^\perp\). Since \(X_2\leq X\leq X\vee Y\), we have
\[
X_2\in \langle X\vee Y\rangle\cap I^\perp=\langle Y\rangle.
\]
Thus \(X_2\leq X\wedge Y=\hat{0}\), and hence \(X_2=\hat{0}\). It follows
that \(X=X_1\in I\), contradicting the assumption that \(X\notin I\).
This completes the proof.
\end{proof}

The proof of Theorem \ref{thm:geometric-realization} proceeds by choosing a translate
\(Y+\alpha\) avoiding all flats in \(L(\A_I)\setminus I\), defining
\(\tau(X)=X\cap(Y+\alpha)\) for \(X\in I\), proving bijectivity, and then
checking that the inverse map is order-preserving.

\begin{proof}[Proof of Theorem \ref{thm:geometric-realization}]
Let \(Y\) be a maximal element of \(I^\perp\) and fix \(y_0\in Y\).
For each \(X\in L(\A_I)\setminus I\), Lemma
\ref{lem:avoidance-for-realization} gives one of two cases. If
\(X\vee Y\) does not exist, equivalently \(X\cap Y=\emptyset\), define
\[
B_X:=X+Y-2y_0
=\{x+y-2y_0:x\in X,\ y\in Y\}.
\]
The direction spaces of \(X\) and \(Y\) cannot span \(V\), since
otherwise \(X\) and \(Y\) would intersect. Hence \(B_X\) is a proper
affine subspace. If instead \(X\wedge Y\neq\hat{0}\), choose a
hyperplane \(H_X\in\A\) containing both \(X\) and \(Y\), and define
\[
B_X:=H_X-y_0.
\]
This is also a proper affine subspace.

In either case,
\begin{equation}\label{eq:bad-translation-containment}
X\cap(Y+\beta)\neq\emptyset\quad\Longrightarrow\quad \beta\in B_X.
\end{equation}
Indeed, write \(x=y+\beta\) with \(x\in X\) and \(y\in Y\). In the
first case, \(2y_0-y\in Y\), so
\(\beta=x+(2y_0-y)-2y_0\in B_X\). In the second, \(x-y\) lies in the
direction space \(H_X-y_0\), so again \(\beta\in B_X\).
Thus every
\begin{equation}\label{eq:generic-translation-set}
\alpha\in\Omega_{I,Y}:=
V\setminus\bigcup_{X\in L(\A_I)\setminus I}B_X
\end{equation}
avoids all flats in \(L(\A_I)\setminus I\). Over an infinite field, a
finite union of proper affine subspaces cannot cover \(V\), so such an
\(\alpha\) exists.

Fix \(\alpha\) satisfying \eqref{eq:generic-translation-set}, and define
$\tau \colon I \longrightarrow L\left((\A_I)^{Y + \alpha}\right)$ by
\[
\tau(X) = X \cap (Y + \alpha).
\]
Since \(I\) is ordered by reverse inclusion, \(\tau\) is
order-preserving. We verify that it is well-defined, surjective,
injective, and order-reflecting.
\begin{enumerate}
    \item[(1)] $\tau$ is well-defined.
For \(X\in I\), the join and rank conditions
\hyperlink{join_condition}{$(\mathrm{I2})$},
\hyperlink{rank_condition}{$(\mathrm{I3})$}, together with
\(X\wedge Y=\hat{0}\), give
\[
r(X\vee Y)=r(X)+r(Y),
\]
or equivalently
\[
d+\dim(X\cap Y)=\dim X+\dim Y.
\]
Thus the direction spaces of \(X\) and \(Y\) span \(V\). Consequently
every translate \(Y+\alpha\) meets \(X\), and the intersection has the
same dimension as \(X\cap Y\):
\[
\dim(X\cap(Y+\alpha))=\dim(X\cap Y).
\]
Moreover, \(X\) is an intersection of hyperplanes in \(\A_I\), and the
same rank argument with \(X=H\in\A_I\) shows that each
\(H\cap(Y+\alpha)\) is a proper hyperplane of \(Y+\alpha\). Hence
\(X\cap(Y+\alpha)\) is a flat of the restricted arrangement, and
\(\tau\) is well-defined.
    \item[(2)] $\tau$ is surjective.
By the definition of the restricted arrangement,
\[
L\left((\A_I)^{Y + \alpha}\right) = \{ X \cap (Y + \alpha) \mid X \in L(\A_I), \; X \cap (Y + \alpha) \neq \emptyset \}.
\]
The choice of \(\alpha\) makes
\(X\cap(Y+\alpha)\) empty for every \(X\in L(\A_I)\setminus I\).
Since \(I\subseteq L(\A_I)\), the preceding display reduces to
\[
L\left((\A_I)^{Y + \alpha}\right) = \{X \cap (Y + \alpha) \mid X \in I,\; X \cap (Y + \alpha) \neq \emptyset \}.
\]
Thus \(\tau\) is surjective.
    \item[(3)]
$\tau$ is injective.
We first show that \(\tau\) is rank-preserving, that is,
\[
r(\tau(X))=r(X)
\]
for all \(X\in I\).
The join and rank conditions give
\[
r(X\vee Y)=r(X)+r(Y),
\]
which in terms of dimensions can be rewritten as
\[
d+\dim(X\cap Y)=\dim(X)+\dim(Y).
\]
\begin{align*}
r(\tau(X))
&= \dim(Y + \alpha) - \dim(X \cap (Y + \alpha)) \\[6pt]
&= \dim(Y) - \dim(X \cap Y)\\[6pt]
&= d - \dim(X) \\[6pt]
&= r(X).
\end{align*}
Next, let \(X_1,X_2\in I\) such that \(\tau(X_1)=\tau(X_2)\), that is,
\[
X_1\cap(Y+\alpha)=X_2\cap(Y+\alpha),
\]
and both are non-empty.
Then
\[
X_1 \cap (Y + \alpha)
= (X_1 \cap X_2) \cap (Y + \alpha)
= X_2 \cap (Y + \alpha).
\]
In particular, \(X_1\cap X_2\) is nonempty and belongs to \(L(\A_I)\).
Since this flat meets
\(Y+\alpha\), the choice of \(\alpha\) forces \(X_1\cap X_2\in I\). Thus
\[
r(\tau(X_1))=r(\tau(X_1\cap X_2))=r(\tau(X_2)).
\]
Since \(\tau\) is rank-preserving, it follows that
\[
r(X_1) = r(X_1 \cap X_2) = r(X_2).
\]
Thus \(X_1=X_1\cap X_2=X_2\), which shows that \(\tau\) is injective.
\item[(4)] $\tau^{-1}$ is order-preserving.
It is enough to show that \(\tau\) reflects the order. Let \(X,Z\in I\)
and suppose that
\[
\tau(X)\leq \tau(Z),
\]
that is, \(Z\cap(Y+\alpha)\subseteq X\cap(Y+\alpha)\) as affine
subspaces. Put \(W=X\cap Z\), equivalently \(W=X\vee Z\) in the
intersection semilattice. Since \(Z\cap(Y+\alpha)\subseteq X\), we have
\[
W\cap(Y+\alpha)=Z\cap(Y+\alpha)\neq\emptyset.
\]
As \(W\in L(\A_I)\), the choice of \(\alpha\) forces \(W\in I\). By the
rank-preserving property already proved,
\[
r(W)=r(\tau(W))=r(\tau(Z))=r(Z).
\]
Moreover \(Z\leq W\). In an intersection semilattice, two comparable flats of
the same rank are equal; hence \(Z=W\), so \(Z\subseteq X\), or
equivalently \(X\leq Z\). Therefore \(\tau^{-1}\) is order-preserving.
\end{enumerate}
We conclude that \(\tau\) is an order-preserving bijection whose inverse
is also order-preserving. Therefore
\[
I \cong L\left((\A_I)^{Y + \alpha}\right)
\]
as required.
\end{proof}

\subsection{The parameter space of translated realizations}

We now make the generic translation set explicit and give a finite-field
criterion.

\begin{theorem}\label{thm:realization-parameter-space}
Let \(\A\) be an essential affine arrangement over a field \(\mathbb{K}\),
let \(I\) be a modular ideal of \(L(\A)\), and let \(Y\) be a maximal
element of \(I^\perp\). For each
\(X\in L(\A_I)\setminus I\), there is a proper affine subspace
\(B_X\subsetneq V\) such that every
\[
\alpha\in\Omega_{I,Y}:=
V\setminus\bigcup_{X\in L(\A_I)\setminus I}B_X
\]
satisfies
\[
I\cong L\bigl((\A_I)^{Y+\alpha}\bigr).
\]
Thus \(\Omega_{I,Y}\) is a Zariski-open subset of \(V\). If
\(\mathbb{K}\) is infinite, it is nonempty.
\end{theorem}

\begin{proof}
The construction of \(B_X\) in the proof of Theorem
\ref{thm:geometric-realization} does not use the infinitude of the field.
Equation \eqref{eq:bad-translation-containment} shows that every
\(\alpha\in\Omega_{I,Y}\) yields the required realization. The set is
Zariski open, and it is nonempty over an infinite field because finitely
many proper affine subspaces cannot cover \(V\).
\end{proof}

Over a finite field this gives the following criterion. Work over
\(\mathbb F_q\) and put
\[
N=|L(\A_I)\setminus I|.
\]
If \(q>N\), then \(\Omega_{I,Y}\neq\emptyset\). Indeed, each bad affine
subspace has at most \(q^{d-1}\) points, so
\[
\left|\bigcup_XB_X\right|\leq Nq^{d-1}<q^d
\]
when \(q>N\). Moreover,
Proposition \ref{prop:graded-modular-ideals} and
\eqref{eq:rank-sum-arrangement} give
\[
\dim(Y+\alpha)=d-r(I^\perp)=r(I).
\]
Thus the constructed realization is essential.

When \(\A\) is central, we may translate a common point to the origin and
assume that \(\A\) is \emph{linear}, meaning that all its hyperplanes pass
through the origin. If \(Z\) is a modular element of \(L(\A)\), Theorem
\ref{thm:geometric-realization} gives a particularly simple description
over an infinite field: every vector outside the hyperplanes containing
\(Z\) is valid; that is, one may take
\[
\alpha\in V\setminus \bigcup_{H\in \A_Z}H.
\]

\begin{coro}\label{cor:central-restriction-realization}
Let $\A$ be a linear arrangement over an infinite field and let $Z$ be a
modular element of $L(\A)$.
Then, for any \(\alpha\in V\setminus\bigcup_{H\in\A_Z}H\),
\[\langle Z \rangle^\perp\cong L(\A^{Z+\alpha}).\]
Consequently,
\[
\chi_{\A}(t) = \chi(\langle Z\rangle,t)
\chi_{\A^{Z+\alpha}}(t).
\]
\end{coro}

\begin{proof}
Propositions \ref{prop:central-modular-ideal} and
\ref{prop:orthogonal-ideal-symmetry} show that
\(\langle Z\rangle\) and \(\langle Z\rangle^\perp\) are modular
ideals.
Apply Theorem \ref{thm:geometric-realization} to
\(\langle Z\rangle^\perp\), choosing \(Z\) as the maximal element
of its orthogonal ideal. This gives
\[
\langle Z \rangle^\perp \cong L\left((\A_{\langle Z \rangle^\perp})^{Z + \alpha}\right)
\]
for every \(\alpha\) satisfying the avoidance condition
\eqref{eq:generic-translation-set}. It remains to show that every
\[
\alpha\in V\setminus\bigcup_{H\in\A_Z}H
\]
satisfies this condition and that
\[
(\A_{\langle Z\rangle^\perp})^{Z+\alpha}
=\A^{Z+\alpha}.
\]
If \(H\in\A_Z\), then \(Z\subseteq H\), and the choice
\(\alpha\notin H\) gives
\((Z+\alpha)\cap H=\emptyset\). Hence such \(H\) does not appear in the
restriction. If \(H\notin \A_Z\), then \(H+Z=V\), so \(H\cap(Z+\alpha)\)
is a nonempty hyperplane of \(Z+\alpha\); equivalently,
\(H\in\A_{\langle Z\rangle^\perp}\).
We next verify the avoidance condition. Fix
\[
X\in
L(\A_{\langle Z\rangle^\perp})
\setminus\langle Z\rangle^\perp.
\]
Lemma \ref{lem:avoidance-for-realization} implies that
\(X\wedge Z\neq\hat{0}\). Hence some \(H\in\A_Z\) contains both \(X\)
and \(Z\). Since \(\alpha\notin H\), the translate \(Z+\alpha\) is
disjoint from \(H\), and therefore from \(X\). Thus
\(X\cap(Z+\alpha)=\emptyset\).
This completes the proof.
\end{proof}

\section{Modular subarrangements and M-ideals}
\label{sec:modular-subarrangements}
Not every modular ideal is induced by a subarrangement. This section
characterizes the induced ones and shows that, after coning, they
correspond to modular flats on the hyperplane at infinity; the same class
is precisely the M-ideals of Bibby and Delucchi \cite{Bibby}.

\begin{prop}\label{prop:subarrangement-principal-intersections}
Let \(I\) be an ideal of \(L(\A)\). Then \(I=L(\A_I)\) if and only if,
for each \(X\in L(\A)\), the ideal \(\langle X\rangle\cap I\) is
principal in \(L(\A)\).
\end{prop}

\begin{proof}
Assume that \(I=L(\A_I)\), and fix \(X\in L(\A)\). To prove that
\(\langle X\rangle\cap I\) is principal, it suffices to show that it has
a unique maximal element. Let \(X_1\) and \(X_2\) be any two maximal
elements of this ideal.
Consider \(X'=X_1\cap X_2\). This intersection is nonempty because
\(X\subseteq X_1\cap X_2\). Since \(X_1,X_2\in L(\A_I)\), their
intersection \(X'\) also lies in \(L(\A_I)=I\). Moreover \(X\subseteq X'\),
so \(X'\in \langle X\rangle\cap I\). The maximality of \(X_1\) and
\(X_2\) gives \(X'=X_1=X_2\). Thus \(\langle X\rangle\cap I\) has a
unique maximal element and is the principal ideal generated by that
element.
	
Conversely, if \(X\in I\), then every hyperplane containing \(X\)
belongs to \(\A_I\), so \(X\in L(\A_I)\). It remains to prove
\(L(\A_I)\subseteq I\).
For any \(Y\in L(\A_I)\), let
\(\langle Y\rangle\cap I=\langle Z\rangle\).
The hyperplanes of \(\A_I\) containing \(Y\) are precisely the
hyperplanes \(H\in\A\) with \(H\leq Y\) and \(H\in I\). Since
\(\langle Y\rangle\cap I=\langle Z\rangle\), these are exactly the
hyperplanes in \(\A_Z\). Therefore
\[
\A_Y \cap \A_I
= \A_Y \cap I
= \A_{\langle Y \rangle \cap I}
= \A_{Z}.
\]
Since \(Y\in L(\A_I)\), it is the intersection, equivalently the join in
the reverse-inclusion order, of all hyperplanes of \(\A_I\) containing it.
Hence
\[
Y = \bigvee_{H \in \A_Y \cap \A_I} H
  = \bigvee_{H \in \A_{Z}} H
  = Z.
\]
It follows that \(Y\in I\). This completes the proof.
\end{proof}

\begin{defi}\label{def:modular-subarrangement}
A subarrangement \(\B\subseteq\A\) is called a \emph{modular
subarrangement} if \(L(\B)\) is a modular ideal of \(L(\A)\).
\end{defi}

Modular subarrangements satisfy the following analogue of
\cite[Lemma 2.27, p.~31]{Orlik-Terao}.

\begin{prop}\label{prop:modular-subarrangement-join-map}
Let \(\B\) be a modular subarrangement of \(\A\). Then for each
\(Y\in L(\B)^\perp\), the map
\[
\tau_Y \colon L(\B) \longrightarrow \{X \vee Y \mid X \in L(\B)\}
\]
defined by \(\tau_Y(X)=X\vee Y\) is an isomorphism of posets.
\end{prop}

\begin{proof}
The map \(\tau_Y\) is order-preserving and surjective by construction.
We prove injectivity and then check that its inverse is order-preserving.
Suppose there exist \(X_1,X_2\in L(\B)\) such that
\(\tau_Y(X_1)=\tau_Y(X_2)\).
Their common image is nonempty, so
\(X_1\cap X_2\cap Y\neq\emptyset\). Hence
\(X_1\vee X_2=X_1\cap X_2\) exists and belongs to \(L(\B)\), and
\begin{align}\label{eq:join-map-injectivity}
X_1 \vee Y = X_2 \vee Y = X_1 \vee X_2 \vee Y.
\end{align}
Since \(L(\B)\) is a modular ideal of \(L(\A)\), the rank condition
\hyperlink{rank_condition}{$(\mathrm{I3})$} applied to
\eqref{eq:join-map-injectivity} gives
\[
r(X_1) + r(Y) = r(X_2) + r(Y) = r(X_1 \vee X_2) + r(Y).
\]
Thus \(r(X_1)=r(X_1\vee X_2)=r(X_2)\). Since both \(X_1\) and \(X_2\)
lie below \(X_1\vee X_2\), equality of ranks gives
\(X_1=X_1\vee X_2=X_2\).
It remains to check that the inverse map is order-preserving. Suppose
\(\tau_Y(X_1)\leq \tau_Y(X_2)\), that is,
\[
X_1\vee Y\leq X_2\vee Y.
\]
Equivalently,
\(X_2\cap Y\subseteq X_1\cap Y\). Thus
\(X_1\cap X_2\cap Y=X_2\cap Y\neq\emptyset\), so
\(X_1\vee X_2\) exists in \(L(\B)\), and
\[
X_1\vee X_2\vee Y=X_2\vee Y.
\]
Applying the rank condition to \(X_1\vee X_2\in L(\B)\) and to
\(X_2\in L(\B)\) gives
\[
r(X_1\vee X_2)+r(Y)=r(X_2)+r(Y),
\]
and hence \(r(X_1\vee X_2)=r(X_2)\). Since \(X_2\leq X_1\vee X_2\),
the two comparable flats have the same rank, so \(X_2=X_1\vee X_2\).
Thus \(X_1\leq X_2\), proving that \(\tau_Y^{-1}\) is order-preserving.
\end{proof}

Let \(\A\) be an affine arrangement in \(\mathbb{K}^d\). Introduce a new
coordinate \(y\). For a hyperplane
\(H\colon \alpha\cdot\boldsymbol{x}=a\) in \(\A\), the \emph{cone}
\(cH\) is the hyperplane in \(\mathbb{K}^{d}\times\mathbb{K}\) given by
\[
cH \colon \alpha \cdot \boldsymbol{x} = a y.
\]
The \emph{cone} \(c\A\) is the linear arrangement
\begin{align}\label{eq:cone-definition}
c\A  = \{cH \mid H\in \A\} \cup \{K\colon y = 0\}.
\end{align}

\begin{fact}\label{fact:cone-poset-isomorphism}
The cone construction induces a rank-preserving poset isomorphism
\[
\theta
 \colon L(\A) \longrightarrow \{Z \in L(c\A) \mid K \nleq Z\},
\]
defined by
\begin{align}\label{eq:cone-poset-isomorphism}
\theta(Y) = \bigvee_{H \leq Y} cH \quad \text{with inverse} \quad \theta^{-1}(Z) = \bigvee_{cH \leq Z} H.
\end{align}
See \cite[Theorem 3.2, p.~371]{Wachs} and
\cite[Exercise (4), p.~428]{Stanley-book}.
\end{fact}

We will use the following standard consequence of the complement
criterion for modular elements \cite[Theorem~1]{Stanley1971}. If \(L\)
is a geometric lattice of rank \(R\) and \(X\in L\) has rank \(a\), then

\begin{equation}\label{eq:complementary-rank-modularity-criterion}
X\text{ is modular}
\quad\Longleftrightarrow\quad
X\wedge Z\neq\hat{0}
\text{ for every \(Z\in L\) of rank \(R-a+1\)}.
\end{equation}

Indeed, a modular \(X\) satisfies
\(r(X\wedge Z)=r(X)+r(Z)-r(X\vee Z)\geq1\). Conversely, if \(X\)
is not modular, Stanley's criterion gives two comparable complements
\(C'<C''\) of \(X\). Semimodularity gives \(r(C')\geq R-a\), so
\(r(C'')\geq R-a+1\). Any element \(Z\leq C''\) of rank \(R-a+1\)
then satisfies \(X\wedge Z=\hat{0}\), contradicting the right-hand
condition in \eqref{eq:complementary-rank-modularity-criterion}.

\begin{theorem}\label{thm:modular-subarrangements-cone}
Let \(\A\) be an affine hyperplane arrangement and \(c\A\) its cone.
Modular subarrangements \(\B\subseteq\A\) are in bijection with modular
elements \(X\in L(c\A)\) satisfying \(K\leq X\). The correspondence is
\[
\B=\varphi(X):=\{H\in\A:cH\leq X\}.
\]
\end{theorem}

\begin{proof}
By Fact \ref{fact:cone-poset-isomorphism}, the map \(\theta\) identifies
\(L(\A)\) with the subposet of \(L(c\A)\) consisting of flats not lying
above \(K\), and it preserves rank.

We first verify that \(\varphi\) is well-defined. Let
\(X\in L(c\A)\) be modular with \(K\leq X\); we show that
\(L(\varphi(X))\) is a modular ideal of \(L(\A)\).
For \(Y\in L(\varphi(X))\) and \(H\in\A\) with \(H\leq Y\),
\eqref{eq:cone-poset-isomorphism} gives
\[
cH \leq \theta(Y) \leq X.
\]
Thus \(H\in\varphi(X)\). If \(Z\leq Y\) in \(L(\A)\), then \(Z\) is
the join of hyperplanes \(H\leq Z\leq Y\), all of which belong to
\(\varphi(X)\). Hence \(Z\in L(\varphi(X))\), proving that
\(L(\varphi(X))\) is an ideal of \(L(\A)\).

We next verify the meet condition. Fix
\(Y\in L(\A)\setminus\{\hat{0}\}\). By Proposition
\ref{prop:subarrangement-principal-intersections}, the ideal
\(\langle Y\rangle\cap L(\varphi(X))\) is principal. If it is nontrivial,
then the meet condition \hyperlink{meet_condition}{$(\mathrm{I1})$}
holds for this \(Y\). If it is trivial, then
\(Y\wedge Z=\hat{0}\) for every \(Z\in L(\varphi(X))\), and hence
\(Y\in L(\varphi(X))^\perp\); in this case
\(\langle Y\rangle\cap L(\varphi(X))^\perp=\langle Y\rangle\) is
nontrivial because \(Y\neq\hat{0}\). Thus the meet condition holds for
\(L(\varphi(X))\).
For \(Y_1\in L(\varphi(X))\) and
\(Y_2\in L(\varphi(X))^\perp\),
\[
K \nleq \theta(Y_1) \leq X \quad \text{and} \quad \theta(Y_2) \wedge X = \hat{0}.
\]
We justify the second assertion. If
\(\theta(Y_2)\wedge X\neq\hat{0}\), then the geometric lattice
\(L(c\A)\) has an atom below \(\theta(Y_2)\wedge X\). This atom cannot
be \(K\), because \(K\nleq\theta(Y_2)\); it is therefore \(cH\) for
some \(H\in\A\).
The relation \(cH\leq X\) gives \(H\in\varphi(X)\), whereas
\(cH\leq\theta(Y_2)\) gives \(H\leq Y_2\). Hence
\(H\in L(\varphi(X))\) and
\(H\wedge Y_2=H\neq\hat{0}\), contradicting
\(Y_2\in L(\varphi(X))^\perp\).
Since \(K\leq X\),
\(\theta(Y_1)\vee K\in\langle X\rangle\). Moreover,
\[
r(\theta(Y_1)\vee K)=r(\theta(Y_1))+1,
\qquad
\theta(Y_2)\in\langle X\rangle^\perp.
\]
By Proposition \ref{prop:central-modular-ideal},
\(\langle X\rangle\) is a modular ideal of \(L(c\A)\).
Therefore
\begin{align*}
r(K \vee \theta(Y_1) \vee \theta(Y_2))
&= r(K \vee \theta(Y_1)) + r(\theta(Y_2))
\tag{by the rank condition \hyperlink{rank_condition}{$(\mathrm{I3})$}} \\[6pt]
&= r(\theta(Y_1)) + r(\theta(Y_2)) + 1
\tag{by \(K\nleq\theta(Y_1)\)} \\[6pt]
&= r(\theta(Y_1) \vee \theta(Y_2)) + 1
\tag{by the rank condition for \(\langle X\rangle\)},
\end{align*}
which implies that \(K\nleq\theta(Y_1)\vee\theta(Y_2)\).
Applying the inverse map in \eqref{eq:cone-poset-isomorphism} gives
\[
Y_1 \vee Y_2 = \theta^{-1}(\theta(Y_1) \vee \theta(Y_2)),
\]
which exists in \(L(\A)\).
Hence the join condition \hyperlink{join_condition}{$(\mathrm{I2})$} holds for \(L(\varphi(X))\).
For the rank condition
\hyperlink{rank_condition}{$(\mathrm{I3})$}, take
\(Y_1\in L(\varphi(X))\) and
\(Y_2\in L(\varphi(X))^\perp\). We have
\(\theta(Y_1)\in\langle X\rangle\) and
\(\theta(Y_2)\in\langle X\rangle^\perp\). The rank condition for
\(\langle X\rangle\), together with the rank-preserving property of
\(\theta\) and the existence of \(Y_1\vee Y_2\) established above,
gives
\[
r(Y_1\vee Y_2)
=r(\theta(Y_1)\vee\theta(Y_2))
=r(Y_1)+r(Y_2).
\]
Thus the rank condition holds for \(L(\varphi(X))\).
Hence \(\varphi(X)\) is a modular subarrangement of \(\A\).
Moreover, every element \(X\) in the domain is recovered from its image by
\[
X=K\vee\bigvee_{H\in\varphi(X)}cH,
\]
and hence \(\varphi\) is injective.

It remains to prove that \(\varphi\) is surjective. Given a modular
subarrangement \(\B\) of \(\A\) of rank \(k\), define
\begin{align}\label{eq:cone-modular-element-from-subarrangement}
X = K \vee \bigvee_{H \in \B} cH.
\end{align}
This is the maximal flat of the cone \(c\B\). The standard rank formula
\(r(c\B)=r(\B)+1\) therefore gives \(r(X)=k+1\).
We first prove that \(X\) is a modular element of \(L(c\A)\).
The lattice \(L(c\A)\) has rank \(d+1\), so
\eqref{eq:complementary-rank-modularity-criterion} says that an element
of rank \(k+1\) is modular if and only if
\[
(c\A)_X\cap(c\A)_Z\neq\emptyset
\]
for every \(Z\in L(c\A)\) of rank \(d-k+1\). Fix such a \(Z\).
\begin{itemize}
    \item If \(K\leq Z\), then \(K\in(c\A)_{X}\cap(c\A)_{Z}\neq\emptyset\).
    \item
If \(K\nleq Z\), then \(Z\) lies in the image of \(\theta\), and
\(r(\theta^{-1}(Z))=d-k+1\). On the other hand,
\(r(L(\B)^\perp)=d-k\), as implied by
\eqref{eq:rank-sum-arrangement}.
Thus,
\begin{align}\label{eq:mideal-rank-inequality}
r(L(\B)^\perp) < r(\theta^{-1}(Z)).
\end{align}
Lemma \ref{lem:join-generation}, applied to the ideal decomposition
\(\{L(\B),L(\B)^\perp\}\), gives
\[
\theta^{-1}(Z)=Y_1\vee Y_2
\]
for some \(Y_1\in L(\B)\) and
\(Y_2\in L(\B)^\perp\). Inequality
\eqref{eq:mideal-rank-inequality} forces \(Y_1\neq\hat{0}\); otherwise
\(\theta^{-1}(Z)=Y_2\) would have rank at most
\(r(L(\B)^\perp)\).
Since \(Y_1\neq\hat{0}\) is a flat of \(L(\B)\), some hyperplane
\(H\in\B\) satisfies \(H\leq Y_1\). Hence
\(cH\leq X\) by \eqref{eq:cone-modular-element-from-subarrangement}. Also \(Y_1\leq \theta^{-1}(Z)\), so the
cone isomorphism \eqref{eq:cone-poset-isomorphism} gives \(cH\leq Z\). Therefore
\(cH \in (c\A)_{X} \cap (c\A)_{Z}\).
\end{itemize}
Thus \((c\A)_X\cap(c\A)_Z\neq\emptyset\) for every \(Z\) of the
required rank. Equivalently, \(X\wedge Z\neq\hat{0}\) for every such
\(Z\), and \eqref{eq:complementary-rank-modularity-criterion} shows that
\(X\) is modular.
It remains to prove that \(\B=\varphi(X)\). The inclusion
\(\B\subseteq\varphi(X)\) follows from the definition of \(X\).
We prove the converse by contradiction:
suppose that there exists \(H\notin\B\) such that \(cH\leq X\) in
\(L(c\A)\).
Then \(H\in L(\B)^\perp\). Indeed, applying the meet condition for
the modular ideal \(L(\B)\) to the rank-one flat \(H\), the alternative
\(\langle H\rangle\cap L(\B)\neq\{V\}\) would force \(H\in L(\B)\), hence
\(H\in\B\), contrary to the choice of \(H\). Therefore
\(\langle H\rangle\cap L(\B)^\perp=\langle H\rangle\), so
\(H\in L(\B)^\perp\).
For any maximal element \(Y\) of \(L(\B)\), Proposition
\ref{prop:graded-modular-ideals} gives \(r(Y)=k\). The join condition for
the modular ideal \(L(\B)\) guarantees that \(Y\vee H\) exists in
\(L(\A)\).
The rank calculation
\begin{align*}
r(\theta(Y) \vee cH \vee K)
&= r(\theta(Y \vee H) \vee K)
\tag{by \eqref{eq:cone-poset-isomorphism}} \\[6pt]
&= r(\theta(Y \vee H)) + 1
\tag{by \(K\nleq\theta(Y\vee H)\)} \\[6pt]
&= r(Y \vee H) + 1
\tag{by \(\theta\)'s rank-preserving property} \\[6pt]
&= k + 2
\tag{by the rank condition \hyperlink{rank_condition}{$(\mathrm{I3})$}}.
\end{align*}
Since \(\theta(Y)\vee cH\vee K\leq X\), this gives
\(r(X)\geq k+2\), contradicting \(r(X)=k+1\). Hence
\(\B=\varphi(X)\).
\end{proof}

\begin{remark*}
Bibby and Delucchi proved that M-ideals of the intersection semilattice
\(L(\A)\) correspond to modular elements of \(L(c\A)\) lying in the
hyperplane \(K\); see \cite[Theorem~4.2.4]{Bibby}. Combining this with
Theorem \ref{thm:modular-subarrangements-cone} gives
\[
\begin{aligned}
\{\text{M-ideals of }L(\A)\}
&\longleftrightarrow
\{\text{modular subarrangements of }\A\}\\
&\subseteq
\{\text{modular ideals of }L(\A)\}.
\end{aligned}
\]
Thus M-ideals correspond precisely to modular ideals induced by
subarrangements; Example \ref{ex:modular-ideal-realization} shows that
the displayed inclusion can be strict.
\end{remark*}

The next proposition gives a restriction on affine-projection
realizations of modular ideals. If the translated realization of \(I\)
is compatible with an affine projection onto its realizing affine
subspace, then \(I\) must come from a subarrangement.

\begin{prop}\label{prop:fibration-implies-mideal}
Let \(I\) be a modular ideal of \(L(\A)\), let \(Y\) be a maximal element
of \(I^\perp\), and choose
\(\alpha\in\Omega_{I,Y}\) as in Theorem
\ref{thm:realization-parameter-space}. Put \(B=Y+\alpha\). Suppose that
there is an affine projection \(\pi\colon V\to B\) such that
\[
H=\pi^{-1}(H\cap B)
\qquad\text{for every }H\in\A_I.
\]
Then \(I=L(\A_I)\). In particular, \(I\) lies in the M-ideal class
described above.
\end{prop}

\begin{proof}
The inclusion \(I\subseteq L(\A_I)\) follows because every flat is the
intersection of the hyperplanes containing it and \(I\) is an ideal.
Conversely, write an arbitrary \(X\in L(\A_I)\) as
\(X=H_1\cap\cdots\cap H_s\), with \(H_i\in\A_I\). The projection
hypothesis gives
\[
X=\pi^{-1}\bigl((H_1\cap B)\cap\cdots\cap(H_s\cap B)\bigr),
\]
so \(X\cap B\neq\emptyset\). By the definition of
\(\Omega_{I,Y}\), every flat in \(L(\A_I)\setminus I\) is disjoint from
\(B\). Hence \(X\in I\), proving \(I=L(\A_I)\). The last assertion
follows from the preceding remark.
\end{proof}

\section{Orlik--Solomon decomposition and subalgebras}
\label{sec:os-factorization}

This section gives the Orlik--Solomon counterpart of the preceding
factorization. We extend Terao's multiplication theorem
\cite[Theorem~3.8]{HT1986} from modular elements of central arrangements
to modular ideals of affine arrangements, and identify exactly when the
two factors are subalgebras.

We first recall the Orlik--Solomon algebra
\cite{OS1980,Orlik-Terao,Falk-OS}. Let \(\A\) be an affine arrangement
and let \(\Bbbk\) be a coefficient field. Let \(E_1(\A)\) be the
\(\Bbbk\)-vector space with basis \(\{e_H\mid H\in\A\}\), and let
\(E(\A)\) be the exterior algebra of \(E_1(\A)\). We write \(uv\) for
\(u\wedge v\); in particular \(e_H^2=0\) and
\(e_He_K=-e_Ke_H\) for \(H,K\in\A\). If \(|\A|=n\), then
\(E(\A)=\bigoplus_{p=0}^n E_p(\A)\), where \(E_0(\A)=\Bbbk\) and
\(E_p(\A)\) is spanned by the monomials \(e_{H_1}\cdots e_{H_p}\).

Define the \(\Bbbk\)-linear map
\(\partial=\partial_E\colon E(\A)\to E(\A)\) by
\(\partial(1)=0\), \(\partial(e_H)=1\), and, for \(p\geq2\),
\[
\partial(e_{H_1} \cdots e_{H_p}) = \sum_{k=1}^{p} (-1)^{k-1} e_{H_1} \cdots \widehat{e}_{H_k} \cdots e_{H_p}
\]
for \(H_1,\ldots,H_p\in\A\). Then \((E(\A),\partial)\) is a chain
complex.

Let \(\mathcal{S}_p(\A)\) be the set of all \(p\)-tuples of hyperplanes
of \(\A\), and set \(\mathcal{S}(\A)=\bigcup_{p\geq0}\mathcal{S}_p(\A)\).
For \(S=(H_1,\ldots,H_p)\in\mathcal{S}(\A)\), write
\[
e_S:=e_{H_1}\cdots e_{H_p},
\qquad
\bigcap S:=H_1\cap\cdots\cap H_p.
\]
For the empty tuple, we use the conventions \(e_S=1\) and
\(\bigcap S=V\).
The tuple \(S\) is \emph{dependent} if \(\bigcap S\neq\emptyset\) and
\(r(\bigcap S)<|S|\), and \emph{independent} if \(\bigcap S\neq\emptyset\) and
\(r(\bigcap S)=|S|\).
Let \(I(\A)\) be the graded ideal of \(E(\A)\) generated by
\[
\{ e_S \mid \bigcap S = \emptyset \} \cup \{ \partial e_S \mid S \text{ is dependent} \}.
\]
The \emph{Orlik--Solomon algebra} of \(\A\) is
\[
A(\A) = E(\A) / I(\A).
\]
For \(X\in L(\A)\), define
\[
\mathcal{S}_X(\A)=\{S\in\mathcal{S}(\A)\mid \bigcap S=X\},
\]
and
\[
E_X(\A)=\sum_{S\in\mathcal{S}_X(\A)} \Bbbk e_S.
\]
Let \(\psi\colon E(\A)\to A(\A)\) be the quotient map and set
\begin{align}\label{eq:os-local-factor}
	A_X(\A) = \psi(E_X(\A)).
\end{align}
Then \(A_X(\A)\) is generated by
\begin{align}\label{eq:os-direct-sum-decomposition}
	\{ \psi(e_S) \mid \text{$S \in \mathcal{S}_X(\A)$ and $S$ is independent} \}.
\end{align}
Every tuple in this generating set has length \(r(X)\); hence
\(A_X(\A)\) is concentrated in degree \(r(X)\).

\subsection{The graded multiplication isomorphism}
\label{subsec:os-linear-factorization}

\begin{lemma}\label{lem:os-induction-tools}
Let \(J\) be a modular ideal of \(L(\A)\).
\begin{enumerate}
\item[{\rm (1)}] Every hyperplane of \(\A\) belongs to exactly one of
\(\A_J\) and \(\A_{J^\perp}\).
\item[{\rm (2)}] If \(R\in\mathcal{S}(\A)\) is dependent, then
\(\psi(e_R)=0\) in \(A(\A)\).
\end{enumerate}
\end{lemma}

\begin{proof}
For (1), apply the meet condition
\hyperlink{meet_condition}{$(\mathrm{I1})$} to a rank-one element
\(H\in L(\A)\). Since \(\langle H\rangle=\{V,H\}\), a nontrivial
principal intersection with either \(J\) or \(J^\perp\) forces
\(H\in J\) or \(H\in J^\perp\). These two possibilities are
exclusive: if \(H\in J\cap J^\perp\), then the definition of
\(J^\perp\), applied to \(H\in J\), would give
\(H\wedge H=V\), impossible for a hyperplane.

For (2), write \(R=(H_0,H_1,\ldots,H_q)\) after reordering if necessary.
If a hyperplane is repeated, then \(e_R=0\) already. We may therefore
assume that the hyperplanes in \(R\) are distinct.
Since \(R\) is dependent, \(\partial e_R\) is one of the generators of
the Orlik--Solomon ideal \(I(\A)\). In the exterior algebra,
\[
e_{H_0}\partial(e_R)=e_R,
\]
because every term except the first contains \(e_{H_0}^2\). Hence
\(e_R\in I(\A)\), and so \(\psi(e_R)=0\).
\end{proof}

\begin{proof}[Proof of Theorem \ref{thm:os-factorization}]
	Set
	\[
	F(\A):=\Big(\bigoplus_{X \in J}A_X(\A)\Big) \otimes
	\Big(\bigoplus_{Y \in J^\perp}A_Y(\A)\Big),
	\]
	the domain of \(\kappa\) in \eqref{eq:os-multiplication-map}.
	
	First observe that \(F(\A)\) and \(A(\A)\) have the same dimension.
	The tensor product \(F(\A)\) is generated by
	\[
	\left\{ \psi(e_S) \otimes \psi(e_T) \mid S \text{ and } T \text{ are independent, } \bigcap S \in J, \text{ and } \bigcap T \in J^\perp \right\}.
	\]
	On these generators,
	\[
	\kappa(\psi(e_S) \otimes \psi(e_T)) = \psi(e_S) \cdot \psi(e_T),
	\]
	so \(\kappa\) is a graded \(\Bbbk\)-linear map. Each \(A_X(\A)\)
	is finite-dimensional, with
	\[
	\dim_{\Bbbk}A_X(\A) = (-1)^{r(X)} \mu(X),
	\]
	and $ A(\A) $ admits a direct sum decomposition
	\begin{align}\label{eq:os-boundary-reduction}
		A(\A) = \bigoplus_{X \in L(\A)} A_X(\A).
	\end{align}
	See \cite[Theorems 3.72 and 3.74, Proposition 3.75; pp.~77--78]{Orlik-Terao}.
	Therefore \(A(\A)\) has dimension
	\[
	\sum_{X \in L(\A)} (-1)^{r(X)} \mu(X) = (-1)^{r(L(\A))}\chi(L(\A), -1).
	\]
	Similarly, \(F(\A)\) has dimension
	\begin{align*}
		\Big( \sum_{X \in J} (-1)^{r(X)} \mu(X) \Big)
		\Big( \sum_{Y \in J^\perp} (-1)^{r(Y)} \mu(Y) \Big),
	\end{align*}
	which equals
	\[
	(-1)^{r(J) + r(J^\perp)} \chi(J, -1) \chi(J^\perp, -1).
	\]
	Since \(\{J,J^\perp\}\) is an ideal decomposition of \(L(\A)\),
	Corollary \ref{cor:arrangement-characteristic-factorization} gives
	\[
	\chi(J, -1) \chi(J^\perp, -1) = \chi(L(\A), -1),
	\]
	and $r(L(\A)) = r(J) + r(J^\perp)$.
	Hence \(F(\A)\) and \(A(\A)\) are finite-dimensional
	\(\Bbbk\)-vector spaces of the same dimension.
	
	It remains to show that \(\kappa\) is surjective.
	By the definition of \(A(\A)\), it is enough to prove
	\(\psi(e_S)\in\operatorname{Im}(\kappa)\) for every independent
	\(S\in\mathcal{S}(\A)\). For the empty tuple this is immediate, since
	\(\psi(e_S)=1=\kappa(1\otimes1)\). We may henceforth assume that
	\(S\) is nonempty.
	
	Fix such a nonempty independent \(S\). Then
	\(\bigcap S\neq\emptyset\), and the assumption that \(S\) is nonempty
	gives \(\bigcap S\neq V=\hat{0}\). The meet condition
	\hyperlink{meet_condition}{$(\mathrm{I1})$} implies that one of
	\(\langle\bigcap S\rangle\cap J\) and
	\(\langle\bigcap S\rangle\cap J^\perp\) is a nontrivial
	principal ideal. By Corollary \ref{cor:orthogonal-involution}, the
	two ideals may be interchanged; this changes products only by the
	usual graded-commutativity sign. Thus, after interchanging \(J\) and
	\(J^\perp\) if necessary, we may assume
	\[
	\langle \bigcap S \rangle \cap J=\langle P\rangle.
	\]
	We shall use Lemma \ref{lem:os-induction-tools} throughout the
	induction below.
	For an independent tuple \(T\in\mathcal{S}(\A)\) with
	\(\bigcap T\leq\bigcap S\), let
	\[
	\ell(T):=\#\{H\text{ occurring in }T:H\in\A_{J^\perp}\}.
	\]
	We prove \(\psi(e_T)\in\operatorname{Im}(\kappa)\) by induction on
	\(\ell(T)\); the desired assertion follows by taking \(T=S\).
	For the base case \(\ell(T)=0\), every hyperplane of \(T\) lies in
	\(\A_J\). Since \(T\) is independent and
	\(\bigcap T\leq\bigcap S\), all hyperplanes of \(T\) lie in
	\(\langle\bigcap S\rangle\cap J=\langle P\rangle\). Their join
	\(\bigcap T\) therefore lies below \(P\), so
	\(\bigcap T\in J\). Hence
	\[
	\psi(e_T) \in A_{\bigcap T}(\A) \subseteq \bigoplus_{X \in J} A_X(\A).
	\]
	Therefore
	\(\psi(e_T)=\kappa(\psi(e_T)\otimes 1)\in\operatorname{Im}(\kappa)\).
	Now assume the assertion holds for all such independent \(T\) with
	\(\ell(T)<m\), and consider \(\ell(T)=m\). We may write
	\(T=(H_1,\ldots,H_p)\), where
	\(H_1,\ldots,H_{p-m}\in\A_J\) and
	\(H_{p-m+1},\ldots,H_p\in\A_{J^\perp}\). Set
	\[
	T_1 = (H_1, \ldots, H_{p-m}) \quad \text{and} \quad T_2 = (H_{p-m+1}, \ldots, H_p).
	\]
	Set \(Z_1=\bigcap T_1\) and \(Z_2=\bigcap T_2\).
	Since every hyperplane of \(T_1\) lies in the principal ideal
	\(\langle\bigcap S\rangle\cap J=\langle P\rangle\), their intersection
	also belongs to \(\langle P\rangle\subseteq J\). Hence \(Z_1\in J\), and
	thus
	\[ \psi(e_{T_1}) \in \bigoplus_{X \in J} A_X(\A). \]
	There are two cases.
	
	(1) If \(Z_2\in J^\perp\), then
	\[
	\psi(e_{T_2})\in\bigoplus_{Y\in J^\perp} A_Y(\A).
	\]
	Therefore,
	\[
	\psi(e_T) = \psi(e_{(T_1, T_2)}) = \psi(e_{T_1}) \cdot \psi(e_{T_2}) = \kappa(\psi(e_{T_1}) \otimes \psi(e_{T_2})) \in \operatorname{Im}(\kappa).
	\]
	
	(2) If \(Z_2\notin J^\perp\), then Lemma \ref{lem:join-generation} gives
	\(Z_2=U\vee W\) for some \(U\in J\) and \(W\in J^\perp\). Since
	\(Z_2\notin J^\perp\), we must have \(U\neq V\). Choose a
	hyperplane \(H_0\in \A\) with \(H_0\leq U\). Since \(J\) is an ideal
	and \(U\in J\), this hyperplane belongs to \(\A_J\), and
	\(H_0\leq U\leq Z_2\).
	The hyperplane \(H_0\) therefore contains \(\bigcap T_2=Z_2\). Adding
	\(H_0\) to \(T_2\) does not change the intersection but increases the
	length by one.
	Moreover, \(H_0\) does not occur in \(T_1\): otherwise the dependent
	tuple \((H_0,T_2)\) would be a subtuple of the independent tuple \(T\),
	which is impossible because every subtuple of an independent tuple is
	independent.
	Hence \( (H_0, T_2) \) is dependent. The same argument, using
	\(\bigcap T\subseteq Z_2\), shows that \( (H_0,T) \) is dependent.
	In \(A(\A)\),
	\[
	0 = \psi(\partial(e_{(H_0, T)})) = \psi\Big( e_T - \sum_{k=1}^{p} (-1)^{k-1} e_{(H_0, H_1, \ldots, \widehat{H_k}, \ldots, H_p)} \Big),
	\]
	so
	\begin{align*}
		\psi(e_T)
		& = \sum_{k=1}^{p} (-1)^{k-1} \psi\Big(e_{(H_0, H_1, \ldots, \widehat{H_k}, \ldots, H_p)}\Big)\\[6pt]
		& = \Big(\sum_{k=1}^{p-m}+\sum_{k=p-m+1}^{p}\Big)
		(-1)^{k-1}\psi\Big(e_{(H_0, H_1, \ldots, \widehat{H_k}, \ldots, H_p)}\Big).
	\end{align*}
	Each term in the first summation is zero in \(A(\A)\), because it
	contains the dependent subtuple \( (H_0, T_2) \). For the second
	summation, any dependent term is zero by Lemma
	\ref{lem:os-induction-tools}(2). It
	therefore remains only to consider an independent tuple
	\[
	(H_0, H_1, \ldots, \widehat{H_k}, \ldots, H_p),
	\qquad p-m+1 \leq k \leq p.
	\]
	This tuple has \(m-1\) hyperplanes in \(\A_{J^\perp}\), and its
	intersection contains \(\bigcap T\), hence is still below \(\bigcap S\) in
	the reverse-inclusion order. The inductive hypothesis implies that
	\[
	\psi\Big(e_{(H_0, H_1, \ldots, \widehat{H_k}, \ldots, H_p)}\Big) \in \operatorname{Im}(\kappa).
	\]
	Therefore, \(\psi(e_T)\in\operatorname{Im}(\kappa)\). Hence
	\(\kappa\) is surjective, which completes the proof.
\end{proof}

\medskip
\noindent The two direct-sum factors in Theorem \ref{thm:os-factorization} are
defined as subspaces of \(A(\A)\). Over an infinite ground field,
Theorem \ref{thm:geometric-realization} realizes \(J\) and
\(J^\perp\) by arrangements \(\B\) and \(\mathcal C\).
Under the rank-preserving flat-poset isomorphisms supplied by that
theorem, corresponding local Brieskorn components have the same
dimensions,
\[
\dim_{\Bbbk}A_X=(-1)^{r(X)}\mu(X).
\]
Taking direct sums therefore identifies the two factors, noncanonically
as graded vector spaces, with \(A(\B)\) and \(A(\mathcal C)\).
Consequently,
\[
A(\A)\cong A(\B)\otimes_{\Bbbk}A(\mathcal C)
\]
as graded vector spaces.

\subsection{Subalgebras and module structure}
\label{subsec:os-algebra-criterion}

For an ideal \(J\subseteq L(\A)\), write
\[
U_J:=\bigoplus_{X\in J}A_X(\A)\subseteq A(\A).
\]
The Brieskorn decomposition implies that if \(u\in A_X(\A)\) and
\(v\in A_Y(\A)\), then \(uv=0\) or
\[
uv\in A_{X\vee Y}(\A),
\]
where the second alternative requires the join to exist. This gives an
exact combinatorial criterion for \(U_J\) to be a subalgebra.

\begin{theorem}\label{thm:subalgebra-criterion}
Let \(J\) be an ideal of \(L(\A)\). The following conditions are
equivalent.
\begin{enumerate}
\item[\rm (1)] \(U_J\) is a graded subalgebra of \(A(\A)\).
\item[\rm (2)] \(J\) is closed under every join that exists in
\(L(\A)\).
\item[\rm (3)] \(J=L(\A_J)\).
\end{enumerate}
If these conditions hold, the inclusion \(\A_J\subseteq\A\) identifies
\(A(\A_J)\) with \(U_J\) as graded algebras.
\end{theorem}

\begin{proof}
Assume (2). The product rule above shows that the product of two
homogeneous local components indexed by \(J\) is either zero or is indexed
by another element of \(J\). Hence \(U_J\) is a graded subalgebra, proving
(2)\(\Rightarrow\)(1).

Assume (1), and take \(X,Y\in J\) for which \(Z=X\vee Y\) exists.
Since \(J\) is an ideal, every hyperplane containing \(X\) or \(Y\)
lies in \(\A_J\). Choose from these hyperplanes an independent tuple
\(S=(H_1,\ldots,H_{r(Z)})\) whose intersection is \(Z\). Such a tuple
exists because the union of the two localized hyperplane sets has
intersection \(X\cap Y=Z\), and therefore has rank \(r(Z)\). Each
\(\psi(e_{H_i})\) belongs to \(U_J\), and closure under multiplication
gives
\[
0\neq\psi(e_S)\in U_J.
\]
The monomial is nonzero because independent monomials are nonzero in the
Orlik--Solomon algebra, while the Brieskorn decomposition places it in
\(A_Z(\A)\). Since the decomposition
\eqref{eq:os-boundary-reduction} is direct, \(U_J\cap A_Z(\A)\neq0\)
forces \(Z\in J\). This proves
(1)\(\Rightarrow\)(2).

If (2) holds, every element of \(L(\A_J)\) is a nonempty intersection,
equivalently an existing join, of hyperplanes in \(\A_J\subseteq J\), and
therefore belongs to \(J\). Conversely, if \(X\in J\), then every
hyperplane \(H\) containing \(X\) satisfies \(H\leq X\) and hence lies in
\(J\). Since \(X\) is the intersection of these hyperplanes,
\(X\in L(\A_J)\). Thus \(J=L(\A_J)\), proving (3). Conversely, if
\(J=L(\A_J)\), any existing join of
elements of \(J\) is an intersection of flats of the subarrangement and
hence lies in \(L(\A_J)=J\). This proves (2)\(\Leftrightarrow\)(3).

Finally, iterating the deletion--restriction exact sequence shows that
the natural map from the Orlik--Solomon algebra of a subarrangement is
injective \cite[Theorem~3.65]{Orlik-Terao}; its image is generated by the degree-one
classes indexed by \(\A_J\). Under (3), its local components are precisely
the components indexed by \(J\), so its image is \(U_J\).
\end{proof}

\begin{coro}\label{cor:modular-subarrangement-free-module}
Let \(J\) be a modular ideal of \(L(\A)\). Then \(U_J\) is a subalgebra
if and only if \(J\) is induced by the subarrangement \(\A_J\), that is,
\[
J=L(\A_J).
\]
In this case \(\A_J\) is a modular subarrangement,
\[
U_J\cong A(\A_J)
\]
as graded algebras, and \(A(\A)\) is a free graded left
\(A(\A_J)\)-module. More precisely, multiplication induces an
\(A(\A_J)\)-module isomorphism
\[
A(\A_J)\otimes_{\Bbbk}U_{J^\perp}
\xrightarrow{\ \sim\ }A(\A).
\]
\end{coro}

\begin{proof}
Theorem \ref{thm:subalgebra-criterion} gives the equivalence and
identifies \(U_J\) with \(A(\A_J)\). Since \(J\) is modular, the equality
\(J=L(\A_J)\) implies that \(\A_J\) is a modular subarrangement by
Definition \ref{def:modular-subarrangement}. The multiplication map of Theorem
\ref{thm:os-factorization} is then \(A(\A_J)\)-linear in the first
tensor factor, and hence is the displayed free-module isomorphism.
\end{proof}

For graded algebras \(R\) and \(S\), write
\(R\widehat\otimes S\) for the graded tensor product, whose multiplication
uses the Koszul sign convention.

\begin{coro}\label{cor:genuine-os-algebra-factorization}
Let \(J\) be a modular ideal. If both \(J\) and
\(J^\perp\) are induced by modular subarrangements, then
multiplication is an isomorphism of graded algebras
\[
A(\A_J)\widehat\otimes A(\A_{J^\perp})
\xrightarrow{\ \sim\ }A(\A).
\]
Conversely, both factors in Theorem \ref{thm:os-factorization} are
subalgebras of \(A(\A)\) exactly when \(J\) and \(J^\perp\) are both
induced by modular subarrangements.
\end{coro}

\begin{proof}
Under the hypotheses, Theorem \ref{thm:subalgebra-criterion} identifies
the two factors with the indicated Orlik--Solomon subalgebras. Their
elements graded-commute inside \(A(\A)\), so multiplication defines a
homomorphism from the graded tensor product. It is bijective by Theorem
\ref{thm:os-factorization}. The converse follows again from Theorem
\ref{thm:subalgebra-criterion}.
\end{proof}

\section*{Declaration on the Use of AI}

During the preparation of this manuscript, the authors used
OpenAI's ChatGPT and Codex to assist with language editing, \LaTeX{}
organization, consistency checks, and the generation and exploration of
candidate mathematical results and refinements.

The earliest version, prepared before this AI-assisted work, is available at
\url{https://arxiv.org/abs/2504.12226v1} and already contained the core
framework and principal results. In the present version, AI-generated
suggestions were used in developing additional results and refinements
concerning simple semimatroids; the Zariski-open parameter space and
finite-field criterion for geometric realizations; fibration and
affine-projection aspects; and subalgebra and graded-algebra results for
the existing Orlik--Solomon decomposition.

All such suggestions and all mathematical content were independently
verified and finalized by the authors, who take full responsibility for the
manuscript.

\end{document}